
\input amstex
\input amsppt.sty
\magnification1200
\vsize=23.5truecm
\hsize=16.5truecm
\vcorrection{-10truemm}
\NoBlackBoxes

\def\Ami{A_{\min}}
\def\Ama{A_{\max}}
\def\inj{\text{\rm i}}
\def\wA{\widetilde A}
\def\pr{\operatorname{pr}}
\def\comega{\overline\Omega }
\def\simto{\overset\sim\to\rightarrow}
\def\ang#1{\langle {#1} \rangle}
\def\crnp{\overline{\Bbb R}^n_+}
\def\rp{ \Bbb R_+}

\document
\topmatter
\title
Krein-like extensions and \\the lower boundedness problem \\for elliptic operators
\endtitle
\author Gerd Grubb \endauthor
\affil
{Copenhagen Univ\. Math\. Dept\.,
Universitetsparken 5, DK-2100 Copenhagen, Denmark.
E-mail {\tt grubb\@math.ku.dk}}\endaffil
\abstract
For selfadjoint extensions $\wA$ of a symmetric densely
defined positive
operator $\Ami$, the lower boundedness problem is the question of
whether 
$\wA$ is lower bounded {\it if and only if} an associated operator $T$
in
abstract boundary spaces is lower bounded. It holds 
when the Friedrichs extension $A_\gamma $ 
has compact inverse (Grubb 1974, also Gorbachuk-Mikhailets 1976); this applies to elliptic operators $A$ on
bounded domains.

For exterior domains, $A_\gamma ^{-1}$ is not compact, and whereas
the lower bounds satisfy $m(T)\ge m(\wA)$, the
implication of lower boundedness from $T$ to $\wA$ has only been known
when $m(T)>-m(A_\gamma )$. We now show it for
general $T$.

The operator $A_a$ corresponding to $T=aI$, generalizing the Krein-von Neumann
extension $A_0$, appears here; its possible lower boundedness for all real  $a$ is
decisive. We study this Krein-like extension, showing for bounded
domains that the discrete eigenvalues satisfy
$N_+(t;A_a)=c_At^{n/2m}+O(t^{(n-1+\varepsilon )/2m})$ for $t\to\infty $.
\comment
For selfadjoint extensions $\wA$ of a symmetric densely
defined positive
operator $\Ami$, the lower boundedness problem is the question of
whether 
$\wA$ is lower bounded {\it if and only if} an associated operator $T$
in
abstract boundary spaces is lower bounded. It holds 
when the Friedrichs extension $A_\gamma $ 
has compact inverse (Grubb 1974, also announced by Gorbachuk and
Mikhailets 1976); this applies to elliptic operators $A$ on
bounded domains.

For exterior domains, $A_\gamma ^{-1}$ is not compact, and whereas
the lower bounds $m(T)$ and $m(\wA)$ generally satisfy $m(T)\ge m(\wA)$, the
implication of lower boundedness from $T$ to $\wA$ has only been known
to hold when $m(T)>-m(A_\gamma )$. We now show it for
general $T$.

The operator $A_a$ corresponding to $T=aI$, generalizing the Krein-von Neumann
extension $A_0$, appears here; its possible lower boundedness for all real  $a$ is
decisive. We study this Krein-like extension for $a\ne 0$, showing for bounded
domains that the discrete eigenvalues satisfy
$N_+(t;A_a)=c_At^{n/2m}+O(t^{(n-1+\varepsilon )/2m})$ for $t\to\infty $; moreover, this holds
for $A_0$.
\endcomment
\endabstract
\endtopmatter
\rightheadtext {Krein-like extensions and lower boundedness}

\subhead 1. Introduction \endsubhead

The study of extensions of a symmetric operator (or a dual pair of
operators) in a Hilbert space has
a long history, with prominent contributions from J.\ von Neumann
in 1929 \cite{N29}, K.\ Friedrichs 1934 \cite{F34}, M.\ G.\ Krein 1947 \cite{K47},
M.\ I.\ Vishik 1952 
\cite{V52}, M.\ S.\ Birman 1956 \cite{B56} and others. The present author made a
number of contributions in 1968--74 \cite{G68}--\cite{G74}, completing the
preceding theories and working out
applications to elliptic boundary value problems, fully for bounded
domains; further developments are found in \cite{G83}, \cite{G84}.

At the same time there was another, separate development of 
abstract extension
theories, where the operator concept gradually began to be 
replaced by the
concept of relations. This development has been
aimed primarily towards applications to ODE, however including 
operator-valued such equations and
Schr\"o{}dinger operators on ${\Bbb R}^n$; keywords in
this connection are:
boundary triples theory, Weyl-Titchmarsh $m$-functions and Krein
resolvent formulas. Cf.\ e.g.\ Ko\v{c}ube\u{\i} \cite{K75}, 
Vainerman \cite{V80},
Lyantze and  Storozh \cite{LS83}, 
Gorbachuk
and Gorbachuk \cite{GG91}, Derkach and Malamud \cite{DM91}, Arlinskii \cite{A99}, Malamud and Mogilevskii \cite{MM02}, 
Br\"uning,
Geyler and Pankrashkin \cite{BGP06}, and their references. In recent years
there have also been applications to elliptic boundary value problems, cf.\ e.g.\ Amrein and Pearson \cite{AP04},
Behrndt and Langer \cite{BL07}, Ryzhov \cite{R07}, Brown, Marletta, Naboko and
Wood \cite{BMNW08}, Gesztesy and Mitrea \cite{GM08}, and their references. 

The connection between the two lines of extension theories has been
clarified in a recent work of Brown, Grubb and Wood
\cite{BGW09}. Further developments for nonsmooth domains are found in
\cite{G08}, Posilicano and Raimondi \cite{PR09}, Gesztesy and Mitrea
\cite{GM11}, Abels, Grubb and Wood \cite{AGW11}.

There still remain some hitherto unsolved questions, for example
concerning operators over exterior (unbounded) sets, and various
questions in spectral theory.

Meanwhile, there have also been developed powerful tools for PDE in microlocal
analysis, beginning with pseudodifferential operators ($\psi $do's)
and, of relevance here,  going on to pseudodifferential boundary operators ($\psi $dbo's) 
with or without parameters. In a modern treatment it is natural to draw on such techniques when they can be applied
efficiently to solve the problems.
Indeed it is the case for the problems treated in the present paper. 
\medskip

\noindent {\it Lower boundedness.} In the study of realizations $\wA$
of a strongly elliptic $2m$-order differential
operator $A$ on a bounded smooth domain $\Omega \subset {\Bbb R}^n$,
it has been 
known since 1974 that the realization is lower bounded if and only if
a certain operator $T$ determining its boundary condition is lower
bounded. (See Grubb \cite{G74}; an announcement for the symmetric case was also 
given by Gorbachuk and 
Mikhailets \cite{GM76}.)
This proof uses the fact that the inverse of the Dirichlet
realization  $A_\gamma
$ (the Friedrichs extension \cite{F34}) is compact. It is a result in functional analysis of operators in Hilbert
space, and in \cite{G74} it is primarily shown in the abstract
setting of closed extensions of dual pairs of lower bounded
operators $\Ami, \Ami'$ with $\Ami\subset \Ama=(\Ami')^*$, as developed in \cite{G68}. Then it is applied to the
study of general normal boundary conditions for strongly elliptic
systems on compact manifolds with boundary.
A further 
analysis of the lower boundedness problem was given in Derkach 
and Malamud \cite{DM91}.

Assuming only positivity of $A_\gamma $, one has
rather easily that lower boundedness of $\wA$ implies lower
boundedness of $T$, and that a conclusion in the opposite direction
holds if the lower bound of $T$ is above minus the lower bound of
$A_\gamma $; the hard question is to treat large negative lower bounds
of $T$.

In the application of the abstract theory to the case where $\Omega $ is an exterior domain (the complement
of a compact smooth set in ${\Bbb R}^n$) the Dirichlet solution operator $A_\gamma
^{-1}$ is not
compact, and it has been an open problem whether one always could conclude
from lower boundedness of $T$ to lower boundedness of $\wA$. We shall
show in this paper that it is indeed so. The proof uses that the
boundary is compact, and takes advantage of principles and
results for 
pseudodifferential boundary operators \cite{B71}, \cite{G84}, \cite{G96}.

Both symmetric and nonsymmetric cases were treated in \cite{G74}, but the decisive
step takes place in the symmetric setting where $\Ami=\Ami'$. Once it is established
there, one can follow the method of \cite{G74} (the passage from Section 2
to Section 3 there) to extend the result to
dual pairs.
Therefore we shall here focus the attention on the symmetric case. 

The abstract theory is recalled in Section 2, its implementation for
exterior domains is explained in Section 3, and the lower boundedness
result is shown in Section 4.

Section 4 ends with some (easier) observations on G\aa{}rding-type 
inequalities, that are not
tied to bounded boundaries in the same way.
\medskip

\noindent {\it Krein-like extensions.} In the treatment of these lower boundedness questions, a certain family of
non-elliptic realizations comes naturally into the picture. They are
generalizations of the Krein-von Neumann extension (\cite{N29}, \cite{K47}) that we shall here
denote $A_0$; it is the restriction of $\Ama$
with domain $D(A_0)= D(\Ami)\dot+ Z$ where $Z=\ker
\Ama$,
and has attracted much interest through the years, see e.g.\ the 
studies of its spectral
properties by Alonso and Simon \cite{AS80}--\cite{AS81}, Grubb \cite{G83}, Ashbaugh,
Gesztesy, Mitrea, Shterenberg and Teschl \cite{AGMST10}, \cite{AGMT10},  with
further references.

The larger family we shall consider (calling them Krein-like
extensions) is the scale of selfadjoint
operators $A_a$ acting as $\Ama$ with domains
$$
D(A_a)=\{u=v+aA_\gamma ^{-1}z+z\mid v\in D(\Ami),\, z\in Z\},\tag1.1
$$
for $a\in{\Bbb R}$. In the application to boundary value problems, they
are determined by Neumann-type boundary
conditions with pseudodifferential elements; however, they are
non-elliptic and the domains contain $L_2$-functions that are not in
$H^s$ for any $s>0$. For both interior and exterior domains,
their lower boundedness is crucial for the general lower
boundedness problem. Moreover, they play a role  in a study \cite{G11} of
perturbations of essential spectra. 

In the case of a bounded domain,
they will have the single point $a$ as essential spectrum, and one can
ask for the asymptotic behavior of the eigenvalue sequence converging
to $+\infty $ that must exist. In the final Section 5, we deal with this
question, showing that the number $N_+(t;A_a)$  of eigenvalues  in
$[r,t]$ (for some $r>a$)
has the asymptotic behavior
$$
N_+(t;A_a)-c_At^{n/2m}=O(t^{(n-1+\varepsilon  )/2m})\text{ for }t\to\infty ,\tag1.2$$
any $\varepsilon >0$, with the same constant $c_A$ as for the
Dirichlet problem. Here we use results for singular Green operators
obtained in \cite{G84}. We also show this estimate for $A_0$.

\subhead 2. The abstract setting  \endsubhead

We first recall how the general
characterization of extensions is set up.

There is given a
symmetric, closed, densely defined
operator $\Ami$ in a complex Hilbert
space $H$, assumed injective with closed range. Moreover, there is given an
invertible selfadjoint extension $A_\gamma $, such that we have
$$A_{\min}\subset A_\gamma \subset A_{\max}\equiv (A_{\min})^*.$$ 
 Let $$\Cal M=\{\wA\mid \Ami\subset\wA \subset A_{\max}\}.$$
To simplify notation, we write $\wA u$ as $Au$, any $\wA\in\Cal M$. Since $\Ami$ has closed
range, there is an orthogonal decomposition
$$
H=R\oplus Z,\quad R=\operatorname{ran} \Ami, \quad Z= \operatorname{ker} \Ama.\tag2.1
$$

When $X$ is a closed subspace of $H$, we denote by $\pr_X u=u_X$
the orthogonal projection of $u$ onto $X$.

The idempotent operators $\pr_\gamma  =A_\gamma ^{-1}\Ama$ and $\pr_\zeta
=I-\pr_\gamma $ on $D(\Ama)$ define a (non-orthogonal) decomposition of $D(\Ama)$
$$
D(\Ama)=D(A_\gamma )\dot+ Z,\tag2.2 
$$
 denoted $u=u_\gamma +u_\zeta
=\pr_\gamma u+\pr_\zeta u$, which 
allows writing an ``abstract Green's formula'' for $u,v\in D(\Ama)$:
$$
(A u,v) - (u, A v)=((A u)_{Z}, v_{\zeta })-(u_\zeta ,
(A v)_Z).\tag2.3
$$
On the basis of (2.3) one can establish
a 1--1 correspondence (\cite{G68}, also described in \cite{G09}, Chapter 13) between the
closed operators $\wA$ in $\Cal M$ and the closed,
densely defined operators  $T\colon V\to W$, where $V$ and $ W$ are closed
subspaces of $Z$, such that
$$
\text{graph of }T=\{ (\pr_\zeta u, (A u)_W ) \mid u\in D(\wA)\}. \tag2.4
$$
 Here $V=\overline{\pr_{\zeta } D(\widetilde A)}$ and 
$W=\overline{\pr_{\zeta } D(\widetilde A^*)}$.
 For a given operator $T\colon V\to W$, one finds the
corresponding operator $\wA$ from the formula
$$
D(\wA)=\{u\in D(\Ama)\mid \pr_\zeta u\in D(T),\;  (Au)_W= T\pr_\zeta u\}.\tag2.5
$$
In this correspondence, one has moreover:
\roster
\item "(a)" $\wA^*$ corresponds analogously to $T^*\colon W\to V$. In
particular, $\wA$ is selfadjoint if and only if $V=W$ and $T=T^*$.
\item "(b)" $\wA$ is symmetric if and only if $V\subset W$ and $T$ is
symmetric. 
\item "(c)" $\operatorname{ker}\wA=\operatorname{ker}T$; 
\quad $\operatorname{ran}\wA=\operatorname{ran}T+(H\ominus W)$.
\item "(d)" When $\wA$ is bijective, 
$$\wA^{-1}=A_\gamma
^{-1}+\inj_{V }T^{-1}\pr_W.\tag2.6 $$
\endroster
Here $\inj_V$ denotes the injection of $V$ into $H$.

The analysis is related to that of Vishik \cite{V52}, except that he sets the
$\wA$ in relation to operators over the nullspace going in the
opposite direction of our $T$'s and in this context focuses on those
$\wA$'s that have closed range.  Our analysis covers all closed $\wA$.

We recall
furthermore that in view of (2.1), the decomposition (2.2) has the refinement 
$$
D(\Ama)=D(\Ami)\dot+ A_\gamma ^{-1}Z\dot+Z;\tag2.7
$$
it allows to show that when $\wA$ corresponds to  $T$, then
$$
D(\wA)=\{u=v+A_\gamma ^{-1}(Tz+f)+z\mid v\in D(\Ami), \, z\in D(T),\,
f\in Z\ominus W \}.\tag2.8
$$

The lower bound of an operator $P$ is denoted by $m(P)$:
$$
m(P)=\inf\{\operatorname{Re}(Pu,u)\mid u\in D(P),\, \|u\|=1\}\ge
-\infty ;\tag2.9 
$$
when it is finite, $P$ is said to be lower bounded.

Assume now moreover that $\Ami$ has a positive lower bound and that
$A_\gamma $ is the Friedrichs extension of $\Ami$; it
has the same lower bound as $\Ami$.
Then we have in addition
the following facts, shown in \cite{G70} (also described in \cite{G09}):
\roster
\item  "(e)" If $m(\wA)>-\infty $, then $V\subset W$ and $m(T)\ge m(\wA)$.
\item  "(f)" If $V\subset W$ and $m(T)> -m(A_\gamma )$, then
$m(\wA)\ge m(T)m(A_\gamma )/(m(T)+m(A_\gamma ))$.
\endroster
The last rule (shown by Birman \cite{B56} for selfadjoint operators $\wA$) is based on the fact that when $V\subset W$, 
$$
(Au,v)=(Au_\gamma ,v_\gamma )+(Tu_\zeta ,v_\zeta ), \text{ for
}u,v\in D(\wA).\tag2.10
$$

The rule (f) does not cover low values of $m(T)$, but this is was overcome in
\cite{G74} when $A_\gamma ^{-1}$ is compact. Here the situation was set in
relation to the situation where the operators are shifted by subtraction
of a spectral parameter $\mu \in \varrho (A_\gamma )$ (the resolvent set), i.e., all realizations $\wA$ are
replaced by $\wA-\mu $. Here we define
$$
Z_\mu =\ker (\Ama-\mu ),\quad \pr_\gamma ^\mu =(A_\gamma -\mu
)^{-1}(\Ama-\mu ),\quad \pr_\zeta ^\mu =I-\pr_\gamma ^\mu,\tag2.11
$$
which gives a decomposition
$$
D(\Ama)=D(A_\gamma )\dot+ Z_\mu \tag2.12 
$$
(note that $D(\wA-\mu )=D(\wA)$, $D(\Ama-\mu )=D(\Ama)$, $D(A_\gamma -\mu )=D(A_\gamma )$). When $\mu $ is real we have, in the same way
as in the case we started out with, a
1--1 correspondence between operators $\wA-\mu $ and
operators $T^\mu \colon V_\mu \to W_\mu $; here $V_\mu =\overline{\pr_{\zeta }^\mu  D(\widetilde A)}$ and 
$W_\mu =\overline{\pr_{\zeta }^\mu  D(\widetilde A^*)}$, and the
properties (a)--(d) have analogues for this correspondence. In
particular, (d) gives a Krein-type resolvent formula when $\mu \in \varrho (\wA)$,
$$(\wA-\mu )^{-1}=(A_\gamma -\mu )
^{-1}+\inj_{V_\mu  }(T^\mu )^{-1}\pr_{W_\mu }; 
$$
there is much more on
this in \cite{BGW09}.

When $\mu <m(A_\gamma )$, $A_\gamma -\mu $ has positive lower
bound $m(A_\gamma )-\mu $,  so also the properties (e) and (f) have analogues in the new
correspondence. In particular, (f) takes the form:
\roster
\item  "(g)" If $V_\mu \subset W_\mu $ and $m(T^\mu )> -(m(A_\gamma )-\mu )$, then
$$m(\wA)-\mu \ge m(T^\mu )(m(A_\gamma )-\mu )/(m(T^\mu )+m(A_\gamma )-\mu ).\tag2.13$$
\endroster
(Here $V\subset W$ implies $V_\mu \subset W_\mu $, see also
Proposition 2.1 below.)
Note the special case:
\roster
\item  "(h)" If $V_\mu \subset W_\mu $ and $m(T^\mu )\ge 0$, then
$m(\wA)\ge \mu $.
\endroster
Hereby the question of whether $\wA$ is lower bounded when $T$ is so, 
is turned into
the question of whether $m(T^\mu )$ becomes $\ge 0$ when
$\mu \to -\infty $.

Define 
$$E^\mu
=\Ama(A_\gamma -\mu )^{-1}=I+\mu (A_\gamma -\mu )^{-1};\tag2.14
$$ it is a homeomorphism in $H$ such that
$$
F^\mu =(\Ama-\mu )A_\gamma ^{-1}=I-\mu A_\gamma ^{-1}\text{ is the
inverse of }E^\mu .\tag2.15
$$
Moreover, $E^\mu $ maps $Z$ homeomorphically onto $Z_\mu $ (with
inverse $F^\mu $). Details are given in \cite{G74} Section 2, where the
following is shown:

\proclaim{Proposition 2.1} Let $\mu <m(A_\gamma )$.
Define the operator $G^\mu $ in $Z$ by
$$
G^\mu =-\mu \pr_Z E^\mu \,\inj_Z,\tag2.16
$$
it is a bounded selfadjoint operator in $Z$.

Let $\wA$ be a closed operator in $\Cal M$, corresponding to $T\colon V\to
W$. Then $\wA-\mu $ corresponds to $T^\mu \colon V_\mu \to W_\mu $,
determined by
$$
\gathered
V_\mu =E^\mu V,\quad W_\mu =E^\mu W,\quad D(T^\mu )= E^\mu D(T).\\
(T^\mu E^\mu v,E^\mu w)=(Tv,w)+(G^\mu v,w)\text{ for }v\in D(T), w\in W.
\endgathered\tag2.17
$$
\endproclaim

Note that in particular, if $V\subset W$,
$$
\operatorname{Re}(T^\mu E^\mu v,E^\mu v)=\operatorname{Re}(Tv,v)+(G^\mu v,v)\text{ for }v\in D(T).
$$
 
One then observes:

\proclaim{Proposition 2.2} The following statements {\rm (i)} and
{\rm (ii)} are equivalent:

{\rm (i)} For any choice of $V\subset W$ and any lower bounded, closed
densely defined operator $T\colon V\to W$ there is a $\mu <m(A_\gamma )$
such that $m(T^\mu )\ge 0$.

{\rm (ii)} For any $t\ge 0$ there is a  $\mu <m(A_\gamma )$
such that $m(G^\mu )\ge t$.

\endproclaim

\demo{Proof} Let (ii) hold, and consider a lower bounded  operator
$T\colon V\to W$; $V\subset W$. Choose $\mu $ such that $m(G^\mu )\ge \max\{-m(T),0\}$.
Then for $v\in D(T)$,
$$
\operatorname{Re}(T^\mu E^\mu v,E^\mu
v)=\operatorname{Re}(Tv,v)+(G^\mu v,v)\ge m(T)\|v\|^2+m(G^\mu
)\|v\|^2\ge 0.
$$
This shows (i).

Conversely, let (i) hold. It holds in particular for the (selfadjoint)
choices
$T=aI$ on $Z$ with $a\in{\Bbb R}$; let $T^\mu _a$ denote the
corresponding operator on $Z_\mu $. By hypothesis there
is a $\mu $ such that $m(T^\mu _a)\ge 0$. Then 
$$
0\le (T_a^\mu E^\mu v,E^\mu
v)=(av,v)+(G^\mu v,v),\tag2.18
$$
and  hence
$$
(G^\mu v,v)\ge -a\|v\|^2, \text{ for all }v\in Z.
$$
To see that (ii) holds for a given $t\ge 0$, we just have to take $a=-t$. \qed
\enddemo

Note that the proof involves the special choice $T=aI$ on $Z$,
corresponding to the Krein-like extension $A_a$, cf.\ (1.1), (2.8). There is
a formulation in terms of those operators, that can immediately be
included:

\proclaim{Proposition 2.3} The two statements {\rm (i)} and {\rm
(ii)} in Proposition {\rm 2.2} are also equivalent with the statement:

{\rm (iii)} For any $a\in{\Bbb R}$, the Krein-like extension $A_a$,
corresponding to the choice $T=aI$ on $Z$, is lower bounded.
\endproclaim

\demo{Proof} The proof of Proposition 2.2 shows that when (i) holds, its
application to the special cases $T=aI$ on $Z$ gives that $m(T^\mu _a)\ge
0$ for $-\mu $ sufficiently large. By the rule (h), $m(A_a)$ then has
lower bound $\ge \mu $. Since $a$ was arbitrary, we conclude that $A_a$ is lower bounded for
any $a\in{\Bbb R}$; hence (iii) holds.

Conversely, when (iii) holds, it assures by the rule (e) applied to
$A_a-\mu $, that for any 
$a$, $m(T^\mu
_a)\ge 0$ for $-\mu $ sufficiently large. This is used in the 
proof of Proposition 2.2 to conclude that $m(G^\mu )$ is then $\ge -a$, implying (ii). \qed
\enddemo

Then \cite{G74} Th.\ 2.12 showed the validity of (i)--(iii) in an important
case:

\proclaim{Theorem 2.4} When $A_\gamma ^{-1}$ is a compact operator in
$H$, then
$$
m(G^\mu )\to\infty \text{ for }\mu \to -\infty .\tag2.19
$$
Consequently, {\rm (i), (ii)} and {\rm (iii)} of Propositions {\rm
2.2} and {\rm 2.3} are valid; and {\rm (f)} can be supplemented with
\roster
\item  "{\rm (f)$'$}" If $V\subset W$ and $m(T)> -\infty $, then
$m(\wA)>-\infty $.
\endroster

\endproclaim 

Also estimates of the type
$$
\operatorname{Re}(Au,u)\ge c\|u\|^2_{\Cal K}-k\|u\|^2_H,\quad u\in
D(\wA),\tag2.20 
$$
were characterized in \cite{G74}, when $D(A_\gamma ^{1/2})\subset \Cal K\subset H$.

The proof of Theorem 2.4 involves a closer study of the 
Krein-like  realizations $A_a$. We return to a
further analysis of them in Section 5. 

We shall now explain how 
the general set-up  is applied to
boundary value problems. Here we focus on exterior problems since
problems for bounded domains were amply treated in \cite{G68}--\cite{G74}.

\subhead 3. The implementation for exterior boundary
value problems \endsubhead

When $\Omega $ is a smooth open subset of ${\Bbb R}^n$ with boundary $\partial\Omega =\Sigma $, we use the standard
$L_2$-Sobolev spaces, with the following notation: $H^s({\Bbb R}^n)$
($s\in 
{\Bbb R}$) has the norm
$\|v\|_s
=\|{\Cal F}^{-1}(\ang\xi ^s{\Cal
F}v)\|_{L_2({\Bbb R}^n)}$; here ${\Cal F}$ is the Fourier transform and
$\ang\xi =(1+|\xi |^2)^{\frac12}$. Next, $H^s(\Omega )=r_{\Omega }H^s({\Bbb
R}^n)$ where $r_\Omega $ restricts to $\Omega $, provided with
the norm $\|u\|_{s}=\inf\{\|v\|_s\mid v\in H^s({\Bbb
R}^n),\, u=r_\Omega v \}$. Moreover,   $H^s_0(\Omega )
=\{u\in H^s({\Bbb R}^n)\mid \operatorname{supp}u\subset \comega \}$;
closed subspace of $H^s({\Bbb R}^n)$. Spaces over the boundary,
$H^s(\Sigma )$, are defined by local coordinates from $H^s({\Bbb
R}^{n-1})$, $s\in{\Bbb R}$. (There are many equally justified
equivalent choices of
norms there; one can choose a particular norm when convenient.) 
When $s>0$, there are dense continuous embeddings
$$
H^{s}(\Sigma )\subset L_2(\Sigma )\subset H^{-s}(\Sigma ),
$$
and we use the customary identification of
$H^{-s}(\Sigma )$ with the antidual space of $H^s(\Sigma )$ (the space
of antilinear, i.e., conjugate linear, functionals), such that the
duality $(\varphi ,\psi )_{-s,s}$ coincides with the $L_2(\Sigma
)$-scalar product when the elements lie there. 
 
Detailed explanations are found in many
books, e.g.\ 
\cite{LM68}, \cite{H63}, \cite{G09}.

In the following, $\Omega $ is primarily considered to be an exterior domain, i.e., the
complement of $\comega_0$, where $\Omega _0$ is a nonempty smooth bounded
subset of ${\Bbb R}^n$. However, the explanations in the following work
equally well for interior domains and for admissible manifolds in the
sense introduced in the book \cite{G96}; this includes smooth domains in
${\Bbb R}^n$ that outside of a
large ball have the form of a halfspace ${\Bbb R}^n_+$
or a cone.

Let $A$ be a symmetric elliptic
operator of order $2m$ on $\Omega $,
$$
Au={\sum}_{|\alpha |,|\beta |\le
m}D^\alpha (a_{\alpha ,\beta }(x)D^\beta u(x)),\quad \overline {a_{\beta ,\alpha }}=a_{\alpha ,\beta },\tag3.1 $$
with complex coefficients $a_{\alpha ,\beta }$ in $C^\infty _b(\comega)$; here $D^\alpha =D_1^{\alpha _1}\cdots
D_n^{\alpha _n}$, $D_j=-i\partial/\partial x_j$, and $C^\infty
_b(\comega )$ denotes the space of $C^\infty $-functions that
are bounded with bounded derivatives of all orders.
The principal symbol $a^0(x,\xi )=
{\sum}_{|\alpha |,|\beta |=
m} a_{\alpha , \beta } \xi ^{\alpha +\beta }  $ is real. $A$ is assumed to
be uniformly strongly elliptic, i.e., $a^0$ satisfies, with $c_1>0$, 
$$
a^0(x,\xi )\ge c_1|\xi |^{2m}, \text{ for }x\in\comega,\,\xi \in
{\Bbb R}^n.\tag3.2  
$$
A typical case of such an operator when $m=1$ is of the form
$$
A=-{\sum}_{j,k=1}^n\partial_ja_{jk}(x)\partial_k+a_0(x)={\sum}_{j,k=1}^nD_ja_{jk}(x)D_k+a_0(x),\tag3.3
$$
with real coefficients satisfying $a_{jk}=a_{kj}$ and
$$
 {\sum}_{j,k}a_{jk}(x)\xi _j\xi _k\ge c_1|\xi |^2,\tag3.4
$$
with $c_1>0$.  

We let $H=L_2(\Omega )$, and as $\Ama$ and $\Ami$ we take the operators acting like $A$ in
$L_2(\Omega )$ and defined by
$$\aligned D(\Ama)&=\{u\in L_2(\Omega )
\mid Au\in L_2(\Omega )\text{ in the distribution sense}\},\\
\Ami&=\text{ the closure of }A|_{C_0^\infty (\Omega )};
\endaligned\tag3.5
$$ 
because of the symmetry,   
$\Ama$ and $\Ami$ are adjoints of one another.
It is well-known (and is accounted for e.g.\ in \cite{G11})
that the strong ellipticity and boundedness estimates imply that the graph-norm $(\|Au\|^2+\|u\|^2)^{\frac12}$ and
the $H^{2m}$-norm are equivalent on $H^{2m}_0(\Omega )$, so 
$$
D(\Ami)= H^{2m}_0(\Omega ).
\tag3.6
$$ 
Moreover, when $A_\gamma $ is taken as the 
Dirichlet realisation of $A$, i.e.,
the restriction of $\Ama$ with
domain $D(\Ama)\cap H^m_0(\Omega )$, then
$$
D(A_\gamma )=H^{2m}(\Omega )\cap H^m_0(\Omega);\tag3.7
$$
and
$A_\gamma $ coincides
with the operator defined by
variational theory (the Lax-Milgram lemma) applied to the sesquilinear
form  with
domain $H^m_0(\Omega )$,
$$
a(u,v)= {\sum}
_{|\alpha |,|\beta |\le m}(a_{\alpha,\beta}D^\beta u, D^\alpha v),\tag3.8
$$
 thus $A_\gamma $ is selfadjoint.
 
We can assume that a large enough constant has been added to $A$ such that 
$$
a(u,u)\ge c_0\|u\|^2, \text{ for }u\in H^m_0(\Omega );\tag3.9
$$
with $c_0>0$;
then $c_0$ is also a lower bound for $\Ami$ and
$A_\gamma $, and $A_\gamma $ is invertible.

The set-up of Section 2 applies readily to these choices of
$\Ami$, $\Ama$ and $A_\gamma $; the operators $\wA\in\Cal M$ are
called realizations of $A$. We shall now recall how the
correspondence between a general $\wA$ and an operator $T\colon V\to W$ is
turned into a charaterization of $\wA$ by a boundary condition.

First we note that there is a Green's formula for $A$, valid for $u,v\in H^{2m}(\Omega )$:
$$
(Au,v)_{L_2(\Omega )}-(u,Av)_{L_2(\Omega )}=(\chi u,\gamma v)_{L_2(\Sigma )^n}-(\gamma u,\chi ' v)_{L_2(\Sigma )^n}.\tag3.10
$$
Here, with $\gamma _ju=(\vec
n\cdot D)^ju|_{\Sigma }$, $\vec n$ denoting the interior normal to the boundary,
$$
\aligned
\gamma u&=\{\gamma _0u,\dots,\gamma
_{m-1}u\},\text{ the Dirichlet data,}\\
  \nu u&=\{\gamma _mu,\dots,\gamma
_{2m-1}u\},\text{ the Neumann data,}\\
\chi u &= \Cal A_{M_0M_1}\nu u+\Cal S\gamma u,\quad \chi 'u=-\Cal A_{M_0M_1}^*\nu u+\Cal S'\gamma u, \text{ Neumann-type data;}
\endaligned\tag3.11
$$
where $\Cal A_{M_0M_1}$
is a certain
skew-triangular invertible matrix of differential operators over
$\Sigma $ derived from $A$, and 
$\Cal S$ and $\Cal S'$
are suitable matrices of differential operators; cf.\ \cite{LM68}, \cite{G71}. In the second-order
case (3.3), one can take $\chi $ and $\chi '$ to be the conormal
derivative $\nu _A$ at the boundary,
$$
\nu _{A}u={\sum}_{j,k}a_{jk}n_j\gamma _0\partial_k u.\tag3.12
$$

Occasionally in the following, we shall use the notation of the
calculus of 
pseudodifferential boundary operators ($\psi $dbo's), as initiated by Boutet de
Monvel \cite{B71} and developed further in e.g.\ \cite{G84}, \cite{G96}; there is
also a detailed introduction in \cite{G09}. The calculus defines Poisson operators $K$ (from
$\Sigma $ to $\Omega $), pseudodifferential trace operators $T$ (from
$\Omega $ to $\Sigma $), singular Green operators $G$ on $\Omega $ (including
operators of the form $KT$) and truncated pseudodifferential operators
on $\Omega $, and their composition rules etc. Since we shall in the 
present paper 
only use final theorems on such
operators, we refrain from taking space up here with a detailed introduction.

Let us introduce the notation
$$
{\Cal
  H}^{s}={\prod}_{0\le j<m}H^{s-j-\frac12}(\Sigma )
,\quad \widetilde {\Cal H}^s={\prod}_{0\le j<m}
H^{s-2m+j+\frac12}(\Sigma );\tag 3.13$$
here $({\Cal
  H}^{s})^*=\widetilde {\Cal H}^{2m-s}$,
$(\widetilde {\Cal H}^{s})^*={\Cal
  H}^{2m-s}$, the dualities denoted
$$
\aligned
&(\varphi ,\psi )_{{\Cal H}^{s\,*}, {\Cal H}^s}\text{ or }
(\varphi  ,\psi )_{\{-s+j+\frac12, s-j-\frac12\}},\\
&(\eta ,\zeta )_{\widetilde{\Cal H}^{s\,*}, \widetilde{\Cal H}^s}\text{ or }
(\eta ,\zeta )_{\{2m-s-j-\frac12, s-2m+j+\frac12\}}.
\endaligned
$$
These dualities are consistent with the scalar product in $L_2(\Sigma
)^n$ when the elements lie there. Note that in particular,
$$
{\Cal H}^0=H^{-\frac12}(\Sigma ),\quad \widetilde{\Cal
H}^0=H^{-\frac32}(\Sigma ),\quad ({\Cal H}^0)^*=H^{\frac12}(\Sigma ),\quad
\text{when }m=1.\tag3.14$$

Denote $D^s_A(\Omega )=\{u\in H^s(\Omega )\mid Au\in
L_2(\Omega )\}$, with norm $(\|u\|^2_s+\|Au\|^2_0)^{\frac12}$. It is
seen as in \cite{LM68} that $C^\infty _{(0)}(\comega)=r_\Omega C_0^\infty ({\Bbb R}^n)$ is dense in
$D^s_A(\comega)$, and it follows from \cite{LM68} that $\gamma $, $\nu $, $\chi $ and $\chi '$ extend to
continuous maps:
$$
\gamma \colon D^s_A(\Omega )\to {\Cal H}^s,\quad \nu  \colon D^s_A(\Omega )\to
{\Cal H}^{s-m},\quad \chi,\chi '  \colon D^s_A(\Omega )\to \widetilde{\Cal
H}^s,\text{ for all }s\in{\Bbb R}. \tag3.15
$$
(The mapping properties are shown in \cite{LM68} for bounded domains, but this
implies (3.15) when the properties are
applied to $\Omega \cap B(0,R)$ for a sphere $B(0,R)$ with
$R$ so large that $\Sigma $ is contained in the interior.)
Moreover, Green's formula continues to hold for these extensions, when $u\in H^{2m}(\Omega )$, $v\in D(\Ama)$: 
$$
(Au,v)-(u,Av)=(\chi u,\gamma
v)_{\{j+\frac12, -j-\frac12\}}-(\gamma u,\chi '
v)_{\{2m-j-\frac12,-2m+ j+\frac12\}}
.\tag3.16
$$
Using that $A_\gamma $ is
invertible, one can moreover show that the nonhomogeneous Dirichlet
problem is uniquely solvable: The mapping
$$
{\Cal A}_\gamma =\pmatrix A\\\gamma \endpmatrix \colon  
H^{s}(\Omega ) \to
\matrix
H^{s-2m}(\Omega )\\ \times \\ {\Cal H}^s\endmatrix\tag3.17
$$
has for $s>m-\frac12$ the
solution operator, continuous in the opposite direction,
$$
{\Cal A}_\gamma^{-1}=\pmatrix R_\gamma& \; K_\gamma\endpmatrix;
 \tag3.18
$$
here $R_\gamma $ is for $s=2m$ the inverse of the Dirichlet 
realization $A_\gamma
$, and $K_\gamma $ is the Poisson operator solving the Dirichlet
problem $Au=0,\gamma u=\varphi $. More documentation is given in \cite{G11}.

Denoting $Z^s_A(\Omega )=\{u\in H^s(\Omega )\mid Au=0\}$ (with
$s$-norm), we have in particular the mapping property for $s>m-\frac12$:
$$
\gamma \colon  Z^s_A(\Omega )\simto {\Cal H}^{s},\tag3.19
$$
it extends to all $s\in {\Bbb R}$. (The extension of
the inverse mapping follows from a general rule for
Poisson operators; the direct mapping is treated as shown
in \cite{LM68}, one may also consult the discussion in \cite{G09}, Chapter 11.)

Denote by $\gamma _Z$ the operator acting like $\gamma $ with precise domain and range
$$
\gamma _Z\colon  Z\simto {\prod}_{ j<m}H^{-j-\frac12}(\Sigma )={\Cal H}^0;\tag3.20
$$
it has an inverse $\gamma _Z^{-1}$ and an adjoint $\gamma _Z^*$ that
map as follows:
$$
\gamma _Z^{-1}\colon {\Cal H}^0\simto Z,\quad \gamma _Z^*\colon ({\Cal
H}^0)^*\simto Z.\tag3.21
$$
Both operators lead to  Poisson operators in the $\psi $dbo
calculus when composed with
$\inj_Z$, 
; here
$\inj_Z\gamma _Z^{-1}$ equals $K_\gamma $.
In the case $m=1$, 
$$
\gamma _Z\colon  Z\simto H^{-\frac12}(\Sigma ),\quad \gamma
_Z^{-1}\colon H^{-\frac12}(\Sigma )\simto Z,\quad \gamma
_Z^*\colon H^{\frac12}(\Sigma )\simto Z.
$$

For the study of general realizations $\wA$ of $A$,
the homeomorphism (3.20) allows us to translate the characterization in
terms of operators $T\colon V\to W$ in Section 2 into a characterization in
terms of operators $L$ over the boundary. 

For $V,W\subset Z$, let $X=\gamma V$, $Y=\gamma W$, with the notation
for the restrictions of $\gamma $:
$$
\gamma _V\colon  V\simto X,\quad \gamma _W\colon  W\simto Y.\tag3.22 
$$
The map $\gamma _V\colon V\simto X$ has the adjoint $\gamma _V^*\colon X^*\simto
V$.
Here $X^*$ denotes the antidual space of $X$, again with a duality
coinciding with the scalar product in $L_2(\Sigma )^n$ when
applied to elements that also belong to 
$ L_2(\Sigma )^n$. The duality is written $(\psi ,\varphi
)_{X^*,X}$. We also write $\overline{(\psi ,\varphi )_{X^*,X}}=(\varphi ,\psi )_{X,X^*}$.
Similar conventions are applied to $Y$.

When $A$ is replaced by $A-\mu $ for $\mu <m(A_\gamma )$, there is a
similar notation where $Z$, $V$ and $W$ are replaced by $Z_\mu $,
$V_\mu $, $W_\mu $. Since $\gamma E^\mu z=\gamma z$ (cf.\
(2.14)),  we have that $\gamma $ defines mappings
$$
\gamma _{V_\mu }\colon V_\mu
\simto X,\quad \gamma _{W_\mu }\colon W_\mu \simto Y, \tag3.23
$$
with {\it the same range spaces $X$ and $Y$ as
when $\mu =0$}. 

We denote $K_{\gamma
  ,X}= \inj_{V }\gamma _V^{-1}\colon X\to V\subset H$, it is a Poisson
operator when $X$ is a product of Sobolev spaces.

Now a given $T\colon V\to W$ is carried over to a closed, densely 
defined operator $L\colon  X\to
Y^*$ by the definition
$$
L=(\gamma _W^{-1})^*T\gamma _V^{-1},\quad D(L)=\gamma _V D(T);\tag3.24
$$
it is expressed in the diagram
$$
\CD
V     @>\sim>  \gamma _V   >    X\\
@VTVV           @VV  L  V\\
   W  @>\sim> (\gamma _W^{-1})^*>   Y^*\endCD 
 \tag3.25$$

Observe that when $v\in D(T)$ and $w\in W$ are carried over to  $\varphi
=\gamma _Vv$ and $\psi =\gamma _Ww$, then $L\varphi
=(\gamma _W^*)^{-1}Tv$ satisfies
$$
(Tv,w)=(L\varphi ,\psi )_{Y^*,Y}.\tag3.26
$$

For the question of semiboundedness we note that when $V\subset W$,
hence $X\subset Y$, then the functionals in $Y^*$ act on the elements of
 $X$. Then when $v\in D(T)\subset V\subset W$, so that $\gamma _Vv=\varphi \in
D(L)\subset X\subset Y$, we may write
$$
(Tv,v)=(L\varphi ,\varphi  )_{Y^*,Y}.\tag3.27
$$
The $L_2$-norm of $v$ is equivalent with the ${\Cal H}^0$-norm of $\varphi $;
$$ 
\|v\|\le c_1\|\varphi \|_{\{-j-\frac12\}}\le c_2\|v\|,\quad \varphi =\gamma _Zv,\tag3.28
$$
for any choice of the equivalent norms (denoted $\|\varphi \|_{{\Cal
H}^0}$ or $\|\varphi \|_{\{-j-\frac12\}}$) on the boundary Sobolev
spaces. (One could also fix the norm, e.g.\ by letting $\gamma _Z$ be an isometry.)
Then 
$$
\operatorname{Re}(Tv,v)\ge c\|v\|^2,\quad v\in D(T),\tag3.29
$$
holds for some $c\in {\Bbb R}$ if and only if 
$$
\operatorname{Re}(L\varphi ,\varphi )_{Y^*,Y}\ge c'\|\varphi
\|^2_{\{-j-\frac12\}},\quad \varphi \in D(L),\tag3.30
$$
holds for some $c'\in {\Bbb R}$, and here $c$ and $ c'$ are
simultaneously $>0$ or $\ge 0$.
(If we fix the norm such that $\gamma _Z$ is an isometry, $c=c'$.) 

The interpretation of the condition in (2.5) as a boundary condition
has been explained in several places, beginning with \cite{G68}, so we can do it
rapidly here. Define the Dirichlet-to-Neumann operator
$$
P^0_{\gamma ,\chi }=\chi \gamma _Z^{-1}=\chi K_\gamma \colon {\Cal
H}^0\to\widetilde {\Cal H}^0;\tag3.31
$$
it is in fact continuous from ${\Cal H}^s$ to $\widetilde {\Cal H}^s$
for all $s\in{\Bbb R}$ because of the mapping properties of $\chi $
and $K_\gamma $. It is a matrix-formed {\it pseudodifferential operator over
$\Sigma $}; this was indicated as plausible in \cite{G68}, and proved in 
detail in \cite{G71} 
on the basis of the work of Seeley on
the Calder\'on projector. It also follows from the general $\psi $dbo calculus.
There is the analogous operator $P^{0}_{\gamma ,\chi '}$, and when the
construction is applied to $A-\mu  $ instead of $A$ we get
the operator
$$
P^\mu  _{\gamma ,\chi }=\chi \gamma _{Z_\mu  }^{-1}.\tag3.32
$$

For $m=1$, 
these operators are of order 1, continuous from
$H^{s-\frac12}(\Sigma )$ to $H^{s-\frac32}(\Sigma )$ for all $s$, and
elliptic of order 1 when $A$ and $\chi $ are chosen as in (3.3),
(3.12).
For higher $m$, the operators are {\it multi-order systems}, of the
form $(P_{jk})_{0\le j,k<m}$ with $P_{jk}$ of order $2m-j-k-1$
(continuous from  $H^{s-k-\frac12}(\Sigma )$ to
$H^{s-2m+j+\frac12}(\Sigma )$ for all $s$). Ellipticity is defined in
relation to the multi-order. When $\Cal S=0$ in (3.11), $P^0 _{\gamma
,\chi }$ is elliptic, meaning that the matrix of principal symbols
$\sigma _{2m-j-k-1}(P_{jk})(x',\xi ')$ is regular for $\xi '\ne
0$. (This follows from the ellipticity of $P^0_{\gamma ,\nu }$ shown
in \cite{G71}, see also \cite{G09}, Ch.\ 11.)

 We now define 
$$
\Gamma ^0=\chi -P^0_{\gamma ,\chi }\gamma ,\quad \Gamma ^{\prime 0}=\chi ' -P^0_{\gamma ,\chi '}\gamma ,\tag3.33
$$
also equal to $\chi A_\gamma ^{-1}\Ama$ resp.\ $\chi 'A_\gamma ^{-1}\Ama$;
they are trace operators in the $\psi $dbo calculus, mapping $D(\Ama)$
(with the graph norm) continuously into $\widetilde{\Cal H}^{2m}=({\Cal H}^0)^*$. They vanish on
$Z$. With these operators there holds a modified Green's formula
$$
(Au,v)-(u,Av)=(\Gamma ^0 u,\gamma
v)_{\{j+\frac12, -j-\frac12\}}-(\gamma u,\Gamma  ^{\prime 0}
v)_{\{-j-\frac12,j+\frac12\}},\tag3.34
$$
valid for all $u,v\in D(\Ama)$. In particular,
$$
(Au,w)=(\Gamma ^0u,\gamma w)_{\{j+\frac12, -j-\frac12\}}, \text{ when }u\in D(\Ama), w\in Z.\tag3.35
$$

When $\wA$ corresponds to $T\colon V\to W$ and $L\colon X\to Y^*$, we can write
$$
(Tu_\zeta ,w)=(T\gamma _V^{-1}\gamma u,\gamma _W^{-1}\gamma
w)=(L\gamma u,\gamma w)_{Y^*,Y},\text{ all }u\in D(\wA), w\in W.\tag3.36
$$
The formula $(Au)_W=Tu_\zeta $ in (2.5) is then turned into
$$
(\Gamma ^0u,\gamma w)_{\{j+\frac12, -j-\frac12\}}=(L\gamma u,\gamma w)_{Y^*,Y},\text{ all }w\in W,
$$
or, since $\gamma $ maps $W$ bijectively onto $Y$,
$$
(\Gamma ^0u,\varphi )_{\{j+\frac12,
-j-\frac12\}}=(L\gamma u,\varphi  )_{Y^*,Y}\text{ for all }\varphi  \in Y.\tag3.37
$$

To simplify the notation, note that the injection $\inj_Y\colon  Y\to {\Cal H}^0$ has as adjoint the mapping
$\inj _{Y}^*\colon  ({\Cal H}^0)^*\to Y^*$ that sends a functional $\psi $
on ${\Cal H}^0$ over into a functional $\inj_Y^*\psi $ on $Y$ by:
$$
(\inj_Y^*\psi ,\varphi )_{Y^*,Y}=(\psi ,\varphi
)_{\{j+\frac12,-j-\frac12\}}\text{ for all }\varphi \in Y.
$$
With this notation (also indicated in \cite{G74} after (5.23)), (3.37) may be rewritten as
$$
\inj_Y^*\Gamma ^0u = L\gamma u,
$$
or, when we use that $\Gamma ^0=\chi -P^0_{\gamma ,\chi }\gamma $,
$$
\inj_Y^*\chi  u= (L+\inj_Y^*P^0_{\gamma ,\chi })\gamma u. \tag3.38
$$

We have then obtained:
 
\proclaim{Theorem 3.1} For a closed operator $\wA\in{\Cal M}$, the
following statements {\rm (i)} and {\rm (ii)} are equivalent:

{\rm (i)} $\wA$ corresponds to $T\colon V\to W$ as in Section {\rm 2}. 

{\rm (ii)} $D(\wA)$ consists of the functions  $u\in D(\Ama)$
that satisfy the boundary condition
$$
\gamma u\in D(L),\quad 
\inj_Y^*\chi  u= (L+\inj_Y^*P^0_{\gamma ,\chi })\gamma u.
\tag3.39
$$
 Here $T\colon V\to W$ and $L\colon X\to Y^*$ are defined from one another as
described in {\rm (3.22)--(3.25)}.

\endproclaim

Note that when
$Y$ is the full space ${\Cal H}^0$, $\inj_{Y^*}$ is superfluous, and (3.39) is a
{\it Neumann-type condition}
$$
\gamma u\in D(L), \quad \chi u=(L+P^0_{\gamma ,\chi })\gamma u.\tag3.40
$$

The whole construction can be carried out with $A$ replaced by $A-\mu
$, when $\mu <m(A_\gamma )$. 
We define $L^\mu $ from $T^\mu $ as
in (3.24)--(3.25) with $T\colon V\to W$ replaced by $T^\mu \colon V_\mu \to W_\mu
$ and use of (3.23); here
$$
L^\mu =(\gamma _{W_\mu }^{-1})^*T^\mu \gamma _{V_\mu }^{-1},\quad D(L^\mu )=\gamma _{V_\mu } D(T)=D(L).\tag3.41
$$
Theorem 3.1 implies:

\proclaim{Corollary 3.2} Let $\mu <m(A_\gamma )$. For a closed
operator $\wA\in{\Cal M}$, the 
following statements {\rm (i)} and {\rm (ii)} are equivalent:

{\rm (i)} $\wA-\mu $ corresponds to $T^\mu \colon V_\mu \to W_\mu $ as in Section {\rm 2}. 

{\rm (ii)} $D(\wA)$ consists of the functions  $u\in D(\Ama)$ such that
$$
\gamma u\in D(L),\quad\inj_Y^*\chi  u= (L^\mu +\inj_Y^*P^\mu _{\gamma ,\chi
})\gamma u. \tag3.42
$$

\endproclaim

Since the boundary conditions (3.39) and (3.42) define the same
realization, we obtain moreover the information that
$$
(L^\mu +\inj_Y^*P^\mu _{\gamma ,\chi
})\gamma u=(L +\inj_Y^*P^0 _{\gamma ,\chi
})\gamma u, \text{ for }\gamma u\in D(L),
$$
i.e.,
$$
L^\mu =L+\inj_Y^*(P^0 _{\gamma ,\chi }-P^\mu  _{\gamma ,\chi })\text{
on }D(L).\tag3.43
$$

\subhead 4. The lower boundedness question  \endsubhead

We have shown in Section 2 that the general conclusion of lower
boundedness from $T$ to $\wA$ (hence from $L$ to $\wA$ in view of
(3.28)--(3.30)) hinges on whether the lower bound of
$G^\mu $ takes arbitrary high values when $\mu \to -\infty $. Let us
identify $G^\mu $ in terms of the operators over $\Sigma $.

\proclaim{Proposition 4.1}  Let $\mu <m(A_\gamma )$. We have that
$$
(G^\mu v,w)=((P^0 _{\gamma ,\chi }-P^\mu _{\gamma ,\chi })\gamma _Z
v,\gamma _Zw)_{\{j+\frac12,-j-\frac12\}},\text{ for }v,w\in Z.\tag4.1
$$
In other words,
$$
G^\mu =(\gamma _Z^*)^{-1}(P^0 _{\gamma ,\chi }-P^\mu _{\gamma ,\chi })\gamma _Z^{-1}.\tag4.2
$$
In particular, $P^0 _{\gamma ,\chi }-P^\mu _{\gamma ,\chi }$ is continuous
from ${\Cal H}^0$ to $({\Cal H}^0)^*=\widetilde{\Cal H}^{2m}$.
\endproclaim

\demo{Proof} This is easily seen by use of  the correspondence between realizations and
operators over the boundary, applied to the Krein-von Neumann
extension:

For the case $T=0$ with $V=W=Z$ (defining the Krein-von
Neumann extension), let us denote the operator corresponding to $A_0-\mu $
in the $\mu $-dependent setting by $T^\mu _0$. Here $L=0$,
continuous from
from ${\Cal H}^0$ to $({\Cal H}^0)^*$, 
and we denote the corresponding $\mu
$-dependent operator by $L^\mu _0$; it is likewise continuous
from ${\Cal H}^0$ to $({\Cal H}^0)^*$. By (3.43),
$$
L^\mu _0=P^0 _{\gamma ,\chi }-P^\mu  _{\gamma ,\chi }\text{ on }{\Cal H}^0.
\tag4.3$$
This shows the asserted continuity. By (2.17),
$$
(G^\mu v,w)=(T_0^\mu E^\mu v,E^\mu w), \text{ for }v,w\in Z.\tag4.4
$$
Then furthermore,
$$
\aligned
(T_0^\mu E^\mu v,E^\mu w)&=(L^\mu _0 \gamma _{Z_\mu }E^\mu v,\gamma
_{Z_\mu }E^\mu w)_{\{j+\frac12, -j-\frac12\}}=(L^\mu _0 \gamma _Zv,\gamma
_Zw)_{\{j+\frac12, -j-\frac12\}}\\
&= ((P^0 _{\gamma ,\chi }-P^\mu  _{\gamma ,\chi })\gamma _Zv,\gamma
_Zw)_{\{j+\frac12, -j-\frac12\}},
\endaligned \tag4.5
$$ 

This shows (4.1), and hence (4.2). 
\qed

\enddemo

\example{Remark 4.2}
This was also observed in \cite{BGW09}, Remark 3.2, formulated in the case
$m=1$, for a nonsymmetric situation with general complex values of $\mu
$ (then adjoints and primed operators enter).

An alternative proof that does not refer to the correspondence between
realizations and operators over the boundary goes as follows:
Recalling that $E^\mu =\Ama(A_\gamma -\mu )^{-1}$ maps $Z$
homeomorphically onto $Z_\mu $, we have
for $v,w\in Z$, $\varphi =\gamma _Zv$, $\psi =\gamma _Zw$,
$$\aligned
&(G^\mu v,w)=-\mu (E^\mu v,w)=-\mu (A(A_\gamma -\mu )^{-1}v,w)\\
&=-\mu (\chi (A_\gamma -\mu )^{-1}v, \gamma _Zw)_{\{j+\frac12,-j-\frac12\}}=-\mu (\chi (A_\gamma -\mu )^{-1}v, \psi )_{\{j+\frac12,-j-\frac12\}},
\endaligned\tag4.6$$
where we have used Green's formula (3.16) and the fact that $\gamma
(A_\gamma -\mu )^{-1}=0$. Now if $v\in H^{2m}(\Omega )$, we can use that $-\mu (A_\gamma -\mu
)^{-1}=I-A(A_\gamma -\mu )^{-1}$ to write
$$\aligned
(G^\mu v,w)&=(\chi (I-A(A_\gamma -\mu )^{-1})v, \psi
)_{\{j+\frac12,-j-\frac12\}}\\
&=(\chi v, \psi
)_{\{j+\frac12,-j-\frac12\}}-
(\chi A(A_\gamma -\mu )^{-1}v, \psi )_{\{j+\frac12,-j-\frac12\}}\\
&=(\chi \gamma _Z^{-1}\varphi , \psi
)_{\{j+\frac12,-j-\frac12\}}-
(\chi E^\mu \gamma _Z^{-1}\varphi , \psi )_{\{j+\frac12,-j-\frac12\}}\\
&=(P^0_{\gamma ,\chi }\varphi , \psi
)_{\{j+\frac12,-j-\frac12\}}-
(\chi \gamma _{Z_\mu }^{-1}\varphi , \psi
)_{\{j+\frac12,-j-\frac12\}}\\
&=((P^0_{\gamma ,\chi }-P^\mu _{\gamma ,\chi })\varphi , \psi
)_{\{j+\frac12,-j-\frac12\}}.
\endaligned
$$
This shows the identity for smooth functions $v$ in the
nullspace. Since the smooth null-solutions are dense in $Z$, the general
statement follows by approximation.

We note in passing that since 
$\chi (A_\gamma -\mu )^{-1}$ is the adjoint of
the $\mu $-dependent Poisson operator $K^\mu _\gamma $ (by Green's
formula), (4.6) also
leads to the alternative formula
$$
P^0_{\gamma ,\chi }-P^\mu _{\gamma ,\chi }=-\mu (K^\mu _\gamma )^*K_\gamma .\tag4.7
$$
\endexample

The question of the behavior of the lower bound of $G^\mu $ is hereby
turned into the question of the lower bound of $P^0 _{\gamma ,\chi
}-P^\mu  _{\gamma ,\chi }$, in relation to the norm on ${\Cal
H}^0$. Note that this difference is a multi-order system of $\psi
$do's where the entries are of order $2m$ lower than the entries in $P^0_{\gamma ,\chi }$.

Now this will be set in relation to a similar family of operators in a
situation where the domain $\Omega $ is replaced by a bounded set.
Choose a large open ball $B(0,R)$ containing ${\Bbb R}^n\setminus\Omega $ in its
interior. Let $\Omega _<=\Omega \cap B(0,R)$; its boundary
 $\Sigma _<$ consists of the two disjoint pieces
$\Sigma $ and  $\Sigma '=\partial
B(0,R)$.
When the whole construction is applied to $A$ on $\Omega _<$,
we get a family of matrix-formed Dirichlet-to-Neumann operators
$ P^\mu _{\gamma , \chi _<}$ on $\Sigma  _<=\Sigma \cup \Sigma '$.

\proclaim{ Proposition 4.3} For the pseudodifferential operators
$P^\mu _{\gamma , \chi _<}$ on $\Sigma _<$,
we have
$$
(( P^0 _{\gamma ,\chi _<}- P^\mu  _{\gamma ,\chi _< })\varphi 
,\varphi )_{\{j+\frac12, -j-\frac12\}}\ge C(\mu
)\|\varphi \|^2_{\{-j-\frac12\}}
$$
for $\varphi \in {\Cal H}^0_<={\prod}_{j<m}H^{-j-\frac12}(\Sigma _<)$,
with 
$$
C(\mu )\to \infty \text{ for }\mu \to -\infty .
$$

\endproclaim

\demo{Proof} This follows from Theorem 2.4, applied to the operator
$G^\mu _<$ defined for this case. The information on the lower bound
of $G^\mu _<$ carries over to the assertion for $P^0 _{\gamma ,\chi
_<}- P^\mu  _{\gamma ,\chi _<}$, since they are related  as
in Proposition 4.1; recall also (3.28)--(3.30). \qed
\enddemo

(This is of course a qualitative statement, which is independent of
how the norm in ${\Cal H}^0$ has been chosen.)

Define
$$
Q^\mu =P^0 _{\gamma ,\chi }-P^\mu _{\gamma ,\chi },\quad 
Q^\mu _1=r_{\Sigma }( P^0 _{\gamma ,\chi _<}- P^\mu  _{\gamma ,\chi _< })e_{\Sigma },\tag4.8
$$
where $e_{\Sigma }$ extends  distributions on $\Sigma $ by 0 on
$\Sigma '$. Since $\Sigma $ and $\Sigma '$ are disjoint closed
manifolds, both $Q^\mu $ and $Q^\mu _1$ are (matrix-formed) $\psi
$do's on $\Sigma $, continuous from  ${\Cal H}^0$ to $\widetilde{\Cal
H}^{2m}$.

\proclaim{Theorem 4.4} The operator norm from  
${\Cal H}^0$ to $\widetilde{\Cal H}^{2m}$ 
of the difference
$$
Q^\mu -Q^\mu _1=P^0 _{\gamma ,\chi }-P^\mu _{\gamma ,\chi }
-r_{\Sigma }( P^0 _{\gamma ,\chi _<}- P^\mu  _{\gamma ,\chi _< })e_{\Sigma }\tag4.9
$$
is bounded for $\mu \to -\infty $.
\endproclaim

\demo{Proof} In this proof we use microlocal details from the
pseudodifferential calculus. Introductions to $\psi $do's can be found in
many textbooks, e.g.\ in \cite{G09}, Chapters 7--8. 

The use of $\psi $do's on the  manifold $\Sigma $ is somewhat
technical, because
they are defined first by Fourier transformation formulas in ${\Bbb
R}^{n-1}$ and then carried over to $\Sigma $ by local coordinates; in
this process there appear a lot of remainder terms that have to be
handled too. 
The heart of our proof lies in the fact that the remainder terms have
much better asymptotic properties than the given operators (are
``negligible''); this is an aspect of the fact that $\psi $do's
are {\it pseudo-local}.

When $P^0_{\gamma ,\chi }$ is constructed from $A$ and the trace
operators, the construction of its symbol takes place in the
neighborhood of  each
point $(x',\xi ')$, $x'\in\Sigma $ (localized) and $\xi ' \in{\Bbb
R}^{n-1}$.
The same holds for $P^0_{\gamma ,\chi _<}$ on $\Sigma $. But in the
localizations at points of $\Sigma $, $A$, $\gamma $ and $\chi $ are
{\it the same} for the two operators, and therefore the resulting complete
symbols of $P^0_{\gamma ,\chi }$ and $P^0_{\gamma ,\chi _<}$ at
a point of $\Sigma $ must be the same, modulo symbols
of order $-\infty $. (This uses that also the constructions in the  
$\psi $dbo calculus are the same for $\Omega $ and $\Omega _<$ at
points of $\Sigma $.)
It follows that 
$$
P^0_{\gamma ,\chi }-r_{\Sigma }P^0_{\gamma ,\chi _<}e_{\Sigma }\text{ is of order }-\infty ,\tag4.10
$$
i.e., the localized symbol of $P^0_{\gamma ,\chi }-r_{\Sigma }P^0_{\gamma
,\chi _<}e_{\Sigma }$ and all its derivatives are $O((1+|\xi '|)^{-N})$
for all $N\in{\Bbb N}$.
Then the operator is bounded as an operator from any $m$-tuple of Sobolev spaces
over $\Sigma $ to any other;
in particular, it is bounded as an operator from ${\Cal H}^0$ to
$\widetilde{\Cal H}^{2m}$.

Now consider the $\mu $-dependent symbols. There is the difficulty here
that the individual operators $P^\mu _{\gamma
,\chi }$ and $P^\mu _{\gamma ,\chi _<}$ have norms that grow with $|\mu
|$ (even as
operators from ${\Cal H}^0$ to $\widetilde {\Cal H}^0$); this is
demonstrated by the simple example of $1-\Delta $ on a half-space
(considered in \cite{G09}, Chapter 9), where $P^\mu _{\gamma ,\chi }$ has symbol 
$-(1+|\xi '|^2+\mu )^{\frac12}$. We shall then use a sharper version of the
device used for $P^0_{\gamma ,\chi }-r_{\Sigma }P^0_{\gamma ,\chi _<}e_{\Sigma }$. Namely that the operators, being
constructed out of the elliptic differential operator $A-\mu $ and the
differential trace operators, have $\mu $-dependent symbols that are
$\psi $do symbols in the $n$ cotangent variables $(\xi ',\eta _n)$
where $\eta _n=|\mu |^{\frac12}$. (This is the ``easy''
parameter-dependent case, said to be {\it of regularity $+\infty $} in
\cite{G96}, {\it strongly polyhomogeneous} in \cite{GS95}.) 

Again, the local constructions of symbols of $P^\mu _{\gamma
,\chi }$ and $P^\mu _{\gamma ,\chi _<}$ have identical ingredients at
the points of $\Sigma $, and we now deduce that the symbols differ by
a symbol in the parameter-dependent class of order $-\infty $, so that
it is $O((1+|\xi '|+|\mu |^{\frac12})^{-N}) $ for all $N\in{\Bbb N}$,
with all its derivatives. Then the symbol and its derivatives are also
$O((1+|\xi '|)^{-N'}(1+|\mu |)^{-N''})$ for all $N',N''\in{\Bbb N}$.
It follows
that 
$$
P^\mu _{\gamma ,\chi }-r_{\Sigma }P^\mu _{\gamma ,\chi _<}e_{\Sigma }\text{ is of order }-\infty ,\text{ with norm }O((1+|\mu |)^{-N}), \text{
any }N,\tag4.11
$$
as an operator from an arbitrary $m$-tuple of Sobolev spaces to another.
In particular, it is bounded as an operator from ${\Cal H}^0$ to
$\widetilde{\Cal H}^{2m}$ with a bound independent of $\mu $.

The assertion on 
$$
Q^\mu -Q^\mu _1=(P^0 _{\gamma ,\chi }-r_{\Sigma }P^0 _{\gamma ,\chi _<}e_{\Sigma })-(P^\mu _{\gamma ,\chi }
-r_{\Sigma }P^\mu _{\gamma ,\chi _<}e_{\Sigma })$$ 
now follows by adding the two parts.  \qed 
\enddemo

We can then conclude:

\proclaim{ Theorem 4.5} In the situation of exterior domains, the pseudodifferential operators
$P^\mu _{\gamma , \chi}$ on $\Sigma $ satisfy
$$
(( P^0 _{\gamma ,\chi }- P^\mu  _{\gamma ,\chi  })\varphi 
,\varphi )_{\{j+\frac12, -j-\frac12\}}\ge C(\mu
)\|\varphi \|^2_{\{-j-\frac12\}}\ \text{ for }\varphi \in {\Cal H}^0,
\tag4.12
$$
for some function $C(\mu )$ satisfying
$$
C(\mu )\to \infty \text{ for }\mu \to -\infty .\tag4.13
$$

It follows that $m(G^\mu )\to\infty $ for $\mu \to -\infty $, and hence:

In the correspondence described in
Theorem {\rm 3.1}, $X\subset Y$ and $L$ is lower bounded, if and only if 
$\wA$ is lower bounded.
\endproclaim

\demo{Proof} Using Proposition 4.3 and Theorem 4.4 we have  for $\varphi \in {\Cal H}^0$
that
$$
\aligned
(( P^0 _{\gamma ,\chi }- P^\mu  _{\gamma ,\chi  })\varphi 
,\varphi )_{\{j+\frac12, -j-\frac12\}}
&=(Q^\mu _1\varphi 
,\varphi )_{\{j+\frac12, -j-\frac12\}}-((Q^\mu _1-Q^\mu )\varphi 
,\varphi )_{\{j+\frac12, -j-\frac12\}}\\
&\ge C(\mu
)\|\varphi \|^2_{\{-j-\frac12\}}-C_1
\|\varphi \|^2_{\{-j-\frac12\}}\ge C'(\mu
)\|\varphi \|^2_{\{-j-\frac12\}};
\endaligned
$$
where $C'(\mu )$ behaves as in (4.13). In view of (4.1)
and (3.28), we conclude that $m(G^\mu )\to\infty $ for $\mu \to
-\infty $. Then the statements (i) and (ii) of Proposition 2.2 are
valid.

Let $\wA$ correspond to 
$T\colon V\to W$ as in the beginning of Section 2, and to $L\colon
X\to Y^*$ as in Theorem 3.1. As noted earlier, $V\subset W$ and $T$ is
lower bounded, if and only if $X\subset Y$ and $L$ is lower bounded. We
have from rule (e) that lower boundedness of $\wA$ implies $V\subset W$ and
lower boundedness of $T$. We can now complete the argument in the
converse direction: When $X\subset Y$ and $L$ is lower bounded, hence
$V\subset W$ and $T$ is lower bounded, then by Proposition 2.2 (i),
there is a $\mu \in{\Bbb R}$ such that $m(T^\mu )\ge 0$, and hence by
rule (h), $m(\wA)\ge \mu $.
\qed

\enddemo

By Proposition 2.3, we have in particular for the Krein-like
extensions:

\proclaim{ Corollary 4.6} In the exterior domain case one has for any
$a\in{\Bbb R}$ that the Krein-like extension $A_a$ defined
by {\rm (1.1)} is lower bounded.
\endproclaim

We recall that this was already known to hold for bounded domains.

\example{Remark 4.7}
The above theorem says nothing about the size of $C(\mu )$. In \cite{G74}
for the interior domain case, we conjectured that $C(\mu )$ may
possibly be shown to be of the order of magnitude $|\mu |^{1/2m}$. 
Calculations on
second-order cases where $A$ has a structure like $D_n^2+B^2$ in product
coordinates near $\Sigma $, confirm that $C(\mu )$ is of the order of 
magnitude $|\mu |^{1/2}$ then. Such calculations might solve the problem also
for domains with unbounded boundary, provided suitable uniform
ellipticity conditions are satisfied. We may possibly return to this in
detail elsewhere.  
\endexample

Let us end this section by some remarks on other lower
boundedness estimates.  It is used in the above proofs that the boundary
$\Sigma $ is compact. There is a more restricted type of lower
boundedness, that can be shown to hold for $\wA$ and  $L$
simultaneously, in uniformly elliptic situations regardless of
compactness of the boundary, 
namely $m$-coerciveness, also known as  the G\aa{}rding inequality.

 Consider a case where $\Omega $ is admissible in the sense of \cite{G96},
as mentioned in the beginning of Section 3. This assures that $\comega
$ is covered by a finite system of local coordinates, some of them for
bounded pieces, some of them for unbounded pieces, carried over to
subsets of ${\Bbb R}^n$ where the part in $\comega$ resp.\
$\partial\Omega $ carries over to bounded resp.\ unbounded subsets of $\crnp$ resp.\ ${\Bbb
R}^{n-1}$, in a controlled way. Detailed explanations are given in
\cite{G96}, including  the still more
general situation of admissible manifolds. All that was
described in Section 2 works in this case; let us also in addition
mention the trace mapping property
$$
\gamma \colon H^r(\Omega )\to {\Cal H}^r\text{ continuously for }r>m-\tfrac12,\tag4.14
$$
and the
interpolation property: When $0<r<m$, there is for any $\varepsilon >0$ a positive
constant $c(\varepsilon )$ such that 
$$
\|u \|^2_{r}\le \varepsilon \|u \|^2_{m}+c(\varepsilon
)\|u \|^2_{0},\text{ for }u \in H^m(\Omega  ).\tag 4.15
$$

\proclaim{Theorem 4.8}  Let $\Omega \subset{\Bbb R}^n$ be an
admissible domain, and let $\wA$ correspond to $T\colon V\to W$ and
$L\colon X\to Y^*$ as in Sections {\rm 2--3}. Then the following
statements (with positive constants $c,c',c''$) are
equivalent:

{\rm (i)} $D(\wA)\subset H^m(\Omega )$ and $\wA$ satisfies a G\aa{}rding inequality
$$
\operatorname{Re}(\wA
  u,u)\ge c\|u\|_m^2-k\|u\|_0^2,\text{ for }u\in D(\wA).\tag 4.16
$$

{\rm (ii)} $D(T)\subset Z\cap H^m(\Omega )=Z^m_A(\Omega )$, $V\subset W$, and  $T$ satisfies a G\aa{}rding inequality
$$
\operatorname{Re}(T
  z,z)\ge c'\|z\|_m^2-k'\|z\|_0^2,\text{ for }z\in D(T).\tag 4.17
$$

{\rm (iii)} $D(L)\subset {\Cal H}^m$,  $X\subset Y$, and  $L$ satisfies a G\aa{}rding inequality
$$
\operatorname{Re}(L\varphi ,\varphi )_{Y^*,Y}\ge c''\|\varphi
\|^2_{\{m-j-\frac12\}} - k''\|\varphi \|^2_{\{-j-\frac12\}}.\tag4.18
$$

\endproclaim

\demo{Proof} This is a straightforward generalization of the proof for the case of bounded domains in \cite{G70},
Prop.\ 2.7, to admissible domains.  

Note first that the statements in (ii) and (iii) are equivalent in view
of (3.27) and the homeomorphisms (3.19). 

Next, we note that (i) implies in particular 
that $\wA$ is lower bounded. Then (i) implies that $V\subset W$ and
hence $X=\gamma V\subset \gamma W=Y$, in view of property (e) in
Section 2. 
Thus (2.10) holds. When (4.16) is valid and $z\in D(T)$, we
can approximate $A_\gamma ^{-1}Tz$ in $m$-norm by a sequence of
functions $v^j\in D(\Ami)$, since $A_\gamma $ is the Friedrichs
extension of $\Ami$. Let $u^j=-v^j+A_\gamma ^{-1}Tz+z$, then $u^j\in
D(\wA)$ in view of (2.8), with $u^j_\gamma =-v^j+A_\gamma ^{-1}Tz$, $u^j_\zeta =z$. Clearly, $u^j\to z$ in
$H^m(\Omega )$ and $u^j_\gamma =-v^j+A_\gamma Tz\to 0$ in $H^m(\Omega
)$. We combine (2.10) with the inequality (4.16) to see that
$$
\operatorname{Re}(Au^j,u^j)=(Au^j_\gamma ,u^j_\gamma )+\operatorname{Re}(Tz,z)
\ge c\|u^j\|^2_m - k\|u^j\|^2_0.
$$
Here the term $(Au^j_\gamma ,u^j_\gamma )$ is equivalent with
$\|u^j_\gamma \|^2_{m}$, so it goes to 0 for $j\to \infty $, so we
conclude that
$$
\operatorname{Re}(Tz,z)\ge c\|z\|_m^2-k\|z\|_0^2.
$$
Thus (i) implies (ii) and hence also (iii).

Now assume that (ii) and (iii) hold. Using (2.10), we find for $u\in
D(\wA)$ that
$$
\aligned
\operatorname{Re}(Au,u)&=(Au_\gamma , u_\gamma
)+\operatorname{Re}(Tu_\zeta, u_\zeta )\\
&\ge c\|u_\gamma
\|_m^2+c'\|u_\zeta \|_m^2-k'\|u_\zeta \|_0^2\ge c''\|u\|_m^2-k'\|u_\zeta \|_0^2,\endaligned\tag4.19 
$$ 
where we have again used that  $(Au_\gamma ,u_\gamma )$ is equivalent with
$\|u_\gamma \|^2_{m}$. To handle the last term, note that choosing $r$ with
$m-\frac12<r<m $, we have that
$$
\aligned
k'\|u_\zeta \|_0^2&\le c_1\|\gamma u_\zeta \|_{\{-j-\frac12\}}^2
= c_1\|\gamma u\|_{\{-j-\frac12\}}^2\le c_2\|\gamma u \|_{\{r-j-\frac12\}}^2 \\
&\le c_3\|u\|_r^2
\le 
\varepsilon c_3\|u\|_m^2+c(\varepsilon )c_3\|u\|_0^2,
\endaligned \tag4.20
$$
where we used (3.19), (4.14) and (4.15). Then (4.19) implies
$$
\operatorname{Re}(Au,u)\ge (c''-\varepsilon c_3)\|u\|_m^2-c(\varepsilon
)c_3\|u \|_0^2,
$$
which shows (i) when $\varepsilon $ is taken sufficiently small.\qed

\enddemo

The papers \cite{G71} and \cite{G74} give a full analysis of the analytical
details required to have (iii) in cases of normal boundary conditions,
for bounded domains and compact manifolds. This involves a condition for
$m$-coerciveness that is a special case of ellipticity of
the boundary condition (the Shapiro-Lopatinskii condition). The
analysis can be extended to
admissible sets with suitable precautions on uniformity of estimates.

We underline that the discussion of lower bounds as in Theorem
4.5 is valid for much more general realizations, and is {\it not} linked
with ellipticity of the boundary condition. An interesting consequence
for questions of spectral asymptotics is that also for nonelliptic
boundary conditions, lower boundedness of $L$ (or $T$) assures that
there is no eigenvalue sequence going to $-\infty $. (For spectral
asymptotics of resolvent differences, see e.g.\  Birman \cite{B62}, 
Birman and Solomyak
\cite{BS80}, Grubb \cite{G84}, \cite{G11}, Malamud \cite{M10}, and
their references.) 

Estimates with other spaces $\Cal K$ in lieu of $H^m(\Omega )$ are
also treated in our early papers.

\subhead 5. Krein-like extensions and their
spectral asymptotics on bounded domains \endsubhead

We here make a closer study of the Krein-like extensions $A_a$ defined
in (1.1), 
corresponding to the choice $T=aI$ in $Z$. 

\proclaim{Proposition 5.1} The realization $A_a$ represents the
boundary condition 
$$
\chi u=C\gamma u,\text{ with } C=a(\gamma _Z^{-1})^*\gamma _Z^{-1}
+P^0_{\gamma ,\chi },\tag 5.1
$$
in the sense that
$$
D(A_a)=\{u\in D(\Ama )\mid \chi u=C\gamma u\}.\tag 5.2
$$

Here 
$(\gamma _Z^{-1})^*\gamma _Z^{-1}$ is a pseudodifferential
operator continuous from ${\Cal H}^s$ to $ \widetilde {\Cal H}^{s+2m}$,
for all $s\in{\Bbb R}$ (and elliptic as such); it is of $2m$ steps lower order than $P^0_{\gamma ,\chi }$.

\endproclaim

\demo{Proof}
We see from (3.24) that $A_a$
corresponds to 
$$
L_a=a(\gamma _Z^{-1})^*\gamma _Z^{-1},\quad D(L)={\Cal H}^0,\tag5.3
$$
so that $A_a$ is defined by the boundary condition in (5.1). 

To
account for the properties of $(\gamma _Z^{-1})^*\gamma _Z^{-1}$ (for
the interested reader), we
use the $\psi $dbo calculus. Note
that 
$(\gamma _Z^{-1})^*\gamma _Z^{-1}$ has the asserted continuity
property for $s=0$, is bijective, and acts
like $(\gamma _Z^{-1})^*\pr_Z\inj_Z\gamma _Z^{-1}$. Here $\inj_Z \gamma _Z^{-1}$ is the Poisson operator $ K_\gamma $,
as noted earlier, and its adjoint $K^*_\gamma =(\gamma _Z^{-1})^*\pr_Z $ is a trace
operator of class 0 in the $\psi $dbo calculus. Then, by the
composition rules,
$$
(\gamma _Z^{-1})^*\gamma _Z^{-1}= K^*_\gamma K_\gamma 
$$
is a
pseudodifferential operator on $\Sigma $; and it has the asserted 
continuity property
for all $s$ since it has it for $s=0$. It is elliptic as an operator
from ${\Cal H}^s$ to $ \widetilde {\Cal H}^{s+2m}$, because it is bijective. \qed
\enddemo

\example{Remark 5.2}
It should be noted that {\it the boundary condition} (5.1) {\it is not 
elliptic} (does not satisfy the appropriate Shapiro-Lopatinskii
condition). In fact, for pseudodifferential Neumann-type boundary conditions
$\chi u=C\gamma u$ it is known that ellipticity holds if and only if
the $\psi $do $L=C-P^0_{\gamma ,\chi }$ is elliptic as an operator
from ${\Cal H}^s$ to $\widetilde {\Cal H}^s$. The actual $L$ equals
$aK_\gamma ^*K_\gamma $, which has principal symbol 0 as an operator
from ${\Cal H}^s$ to $\widetilde {\Cal H}^s$, since it is of lower
order.

For $m=1$, $C$ is of order 1, continuous from $H^{s-\frac12}(\Sigma )$
to $H^{s-\frac32}(\Sigma )$, and $L=aK_\gamma ^*K_\gamma $ is of order
$-1$, continuous from $H^{s-\frac12}(\Sigma )$ to
$H^{s+\frac12}(\Sigma )$, for all $s$.
\endexample

We henceforth take $a\in
{\Bbb R}\setminus\{0\}$. From (2.6) we then have 
$$
A_a^{-1}=A_\gamma ^{-1}+a^{-1}\pr_Z.\tag5.4
$$
(We here read $\pr_X$ as a mapping in $H$ instead of as a mapping from $H$
to $X$; this will often be the case in the following, and the
meaning should be clear from the context.)

Let us assume from now on, instead of the primary hypothesis for Sections 3--4,
that $\Omega $ is a bounded smooth subset of ${\Bbb R}^n$ with
boundary $\Sigma $; aside from this we keep the notation. As remarked in
the beginning of Section 3, the
explanations there hold also for this case (are in fact easier
to verify).

Since the embedding of $D(A_\gamma )=H^{2m}(\Omega )\cap H^m_0(\Omega )$ into $L_2(\Omega )$ is compact,
the inverse $A_\gamma ^{-1}$ is a compact operator in $L_2(\Omega )$, so $A_\gamma $ has a discrete
spectrum consisting of eigenvalues going to $\infty $. It is well-known
(cf.\ e.g.\ H\"ormander \cite{H85}, Ch.\ 29.3), that 
the counting function $N(t;A_\gamma )$, counting the number of
eigenvalues of $A_\gamma $ in $[0,t]$ with multiplicities, has the asymptotic
behavior 
$$
N(t;A_\gamma )-c_At^{n/2m}=O(t^{(n-1)/2m})\text{ for }t\to \infty;\tag5.5
$$
 here 
$$
c_A=(2\pi )^{-n}\int_{x\in \Omega ,\,a^0(x,\xi )<1}\,dxd\xi .\tag5.6
$$
Equivalently, the $j$'th eigenvalue $\mu _j(A_\gamma ^{-1})$ of
$A_\gamma ^{-1}$ satisfies  
$$
\mu _j(A_\gamma ^{-1})-c'_Aj^{-2m/n}=O(j^{-(2m+1)/n})\text{
for }j\to\infty ;\text{ with }c'_A=c_A^{2m/n}.\tag5.7
$$
(The passage between counting function
estimates and eigenvalue estimates is recalled below in Lemma 5.4 and
its corollary.)

Since $Z$ is infinite dimensional, $a^{-1}\pr_Z$ has the point $a^{-1}$ as
essential spectrum, so $A_a^{-1}$ has essential spectrum consisting of
the points $a^{-1}$ and 0, and $A_a$ has the essential spectrum $\{a\}$.
Since $A_a$ is selfadjoint and not upper bounded (since it extends
$\Ami$), there
must be a sequence of discrete eigenvalues (with finite dimensional
eigenspaces)  above $a$ going to $\infty $. We shall investigate this sequence.

The Krein-von Neumann extension $A_0$ has essential spectrum
$\{0\}$ and an eigenvalue sequence going to $\infty $, and the question of the asymptotic behavior
of that sequence was raised in Alonso and Simon \cite{AS80} and answered in
Grubb \cite{G83}. The result 
was a rather precise estimate of
the function $N_+(t;A_0)$ 
counting the number of eigenvalues in $\,]0,t]$:
$$
N_+(t;A_0)-c_At^{n/2m}=O(t^{(n-\theta )/2m})\text{ for }t\to \infty;\tag5.8 
$$
here $c_A$ is the same constant as for the Dirichlet problem and 
$$
\theta =\max\{\tfrac12-\varepsilon , 2m/(2m-n+1)\}.\tag5.9
$$
We note in passing that the value $\tfrac12-\varepsilon $ came from
the application of an estimate
announced by Kozlov in \cite{K79}, whereas his later paper \cite{K83}, not
available to the author when \cite{G83} was written, has the value
$\frac12$, so (5.9) can immediately be replaced by
$$
\theta =\max\{\tfrac12 , 2m/(2m-n+1)\}.\tag5.10
$$
We show at the end of this section that  the estimate can be improved even
further, to $\theta =1-\varepsilon $ (following up on a remark at the
end of \cite{G83}). This comes after our deduction of a similar estimate
for the operators $A_a$, $a\ne 0$.

The proof of (5.8) was based on a transformation of the eigenvalue equation
$$
A_0u=\lambda u,\text{ with }\lambda \ne 0, \, u\ne 0,\tag5.11
$$
into the problem for the $4m$-order operator $A^2$:
$$
A^2v=\lambda Av\text{ for }v\in H^{2m}_0(\Omega ),\tag5.12
$$
where $u$ and $v$ are  recovered from one another by
$$
v=A_\gamma ^{-1}Au,\quad u=\tfrac1\lambda Av.\tag5.13
$$
There were earlier eigenvalue estimates 
for implicit eigenvalue problems as in (5.12) (as initiated by Pleijel \cite{P61},
surveyed in Birman and Solomyak \cite{BS77}) giving the principal
asymptotics, and the
sharper estimates in (5.8) were obtained by turning the problem into
the study of eigenvalues of the compact operator 
$$
S_0=R_\varrho ^{1/2}\,A \,R_\varrho ^{1/2},\tag5.14
$$
where $R_\varrho $ is the solution operator for the Dirichlet problem
for $A^2$. (Further developments of the implicit eigenvalue problem
are described in \cite{G96}, Ch.\ 4.6.)

The study of $A_0$ has been taken up again recently by Ashbaugh,
Gesztesy, Mitrea, Shterenberg and Teschl \cite{AGMST10}, \cite{AGMT10}, also for nonsmooth domains, with much
additional information. In particular they observe that when
$A=-\Delta $, (5.12) is
of interest as the ``buckling problem'' in elasticity. 

Unfortunately, in the case of $A_a$, we do not have an equally simple
reduction of the eigenvalue problem. Let $u=v+aA_\gamma ^{-1}z+z$ as
in (1.1); then applications of powers of $A$ give 
$$
\aligned
&Au-\lambda u=Av+az-\lambda (v+aA_\gamma ^{-1}z+z)=(A-\lambda
)v+(a-\lambda -a\lambda A_\gamma ^{-1})z,\\
&A^2u-\lambda Au=A^2v-\lambda (Av+az)=(A^2-\lambda A)v-a\lambda z,\\
&A^3u-\lambda A^2u=A^3v-\lambda A^2v.
\endaligned
\tag5.15
$$
We see from the third line that in order for $u$ to be an eigenvector,
$v$ must be an eigenvector of a certain implicit problem for $A^3$. 
Here $A^3$ is
of order $6m$, and the information $v\in H^{2m}_0(\Omega )$ does not give
enough boundary conditions to define an elliptic realization of
$A^3$. But there is a supplementing boundary condition depending on
$\lambda $:

\proclaim{Theorem 5.3} Let $u\in D(A_a)$, with $u=v+aA_\gamma
^{-1}z+z$, $v\in H^{2m}_0(\Omega )$, $z\in Z$. Then $u$ is a nonzero
eigenfunction for $A_a$ with eigenvalue $\lambda \ne a$ if and only if
$v$ is a nonzero solution of the elliptic problem
$$
A^3v=\lambda A^2v,\quad
\gamma v=\nu v=0,\quad \gamma A^2v=\lambda ^2(\lambda -a)^{-1}\gamma Av,
\tag5.16
$$
and 
$$
 z=K_\gamma (\lambda -a)^{-1}\gamma Av.\tag5.17
$$

In particular, $u$, $v$ and $z$ are in $C^\infty (\comega)$ then.
\endproclaim
\demo{Proof}
Assume that $Au=\lambda u$, $\lambda \ne a$. It follows from (5.15)
that then $A^3v=\lambda A^2v$. Since $v\in H^{2m}_0(\Omega )$,  $\gamma v=\nu v=0$ (recall
(3.11)). From the first line in (5.15) it is seen that
$$
Av=\lambda (v+aA_\gamma ^{-1}z+z)-az,
$$
which implies
$$
\gamma Av=(\lambda -a)\gamma z,\text{ hence }\gamma z=(\lambda -a)^{-1}\gamma Av.\tag5.18
$$
Moreover,
$$
A^2v=A(\lambda v+\lambda aA_\gamma ^{-1}z+(\lambda -a)z)=\lambda Av+\lambda az,
$$
and hence
$$
\gamma A^2v=\lambda \gamma Av+\lambda a\gamma z=(\lambda +\lambda a(\lambda -a)^{-1})\gamma Av=\lambda ^2(\lambda -a)^{-1}\gamma Av.
$$
This shows the last boundary condition in (5.16) for $v$.
We also see from (5.18) that  $z$ is determined from $v$ by
$z=K_\gamma (\lambda -a)^{-1}\gamma Av$, showing (5.17). Clearly $u\ne
0$ implies $v\ne 0$.

Conversely let $v$ be a nontrivial solution of (5.16), define $z$ by
(5.17) and let $u=v+aA_\gamma ^{-1}z+z$. By the third line of (5.15),
the function $f=A^2u-\lambda Au$ satisfies
$Af=0$; moreover, by the second line,
$$
\gamma f=\gamma (A^2v-\lambda Av-a\lambda z)=\gamma A^2v-\lambda \gamma Av-a\lambda (\lambda -a)^{-1}\gamma Av=0;
$$
where we used (5.17) and the last boundary condition in (5.16). Then
by the unique solvability of the Dirichlet problem, $f=0$. 

Now let $g=Au-\lambda u$, then $Ag=f=0$, and, by the first line of (5.15), 
$$
\gamma g=\gamma (A-\lambda )v+\gamma (a-\lambda )z=\gamma
Av+(a-\lambda )(\lambda -a)^{-1}\gamma Av=0, 
$$
so $g=0$. This shows that $Au=\lambda u$.

The problem is elliptic, since it is a perturbation by lower order
terms of the problem
$$
A^3v=0,\quad \gamma v=\nu v=0,\quad \gamma A^2v=0,
$$
 which only has the zero solution (indeed, $A^3v=0$ and $\gamma
A^2v=0$ imply $A^2v=0$, and then $\gamma v=\nu v=0$ implies $v=0$).
Then since there are $3m$ boundary conditions of different orders, the problem is
elliptic. In particular, the solution of (5.16) is in $C^\infty (\comega)$. \qed 
\enddemo

There may possibly be a strategy to find spectral asymptotics formulas
for the very implicit eigenvalues $\lambda $ of (5.16). 
But rather than pursuing this, we shall apply functional analytical
methods to $A_a$ combined with  $\psi $dbo results, using perturbation
theory for the identity (5.4).

Let us first show how the asymptotic behavior of the counting
functions for positive eigenvalues is related to the asymptotic behavior of 
positive eigenvalues of the inverse operator.

\proclaim{Lemma 5.4} Let $P$ be a selfadjoint invertible operator
whose spectrum on $\rp$ is discrete, consisting of a nondecreasing sequence of positive eigenvalues
$\lambda _{j,+}(P)$ going to  $\infty $ for $j\to\infty $ (repeated
according to multiplicities). Let $N_+(t;P)$ denote
the number of eigenvalues in $[0,t]$, and let $\mu _{j,+}(P^{-1})=\lambda
_{j,+}(P)^{-1}$. Let $C>0$ and let $\beta >\alpha >0$.

There exists $c_1>0$ such that 
$$
|\mu _{j,+}(P^{-1})-Cj^{-\alpha }|\le c_1j^{-\beta }\text{ for all
}j\in{\Bbb N},\tag5.19 
$$
if and only if there exists $c_2>0$ such that
$$
|N_+(t;P)-C^{1/\alpha }t^{1/\alpha }|\le c_2t^{(1+\alpha -\beta )/\alpha
}\text{ for all }t>0.\tag5.20
$$
 
\endproclaim

\demo{Proof} This goes as in the proof for the compact case in \cite{G78}, Lemma 6.2 (a very detailed version is given in \cite{G96}, Lemma A.5): Rewrite (5.19) as
$$
|C^{-1}j^\alpha \mu _{j,+}(P^{-1})-1|\le c_3j^{\alpha -\beta },
$$
$c_3=c_1C^{-1}$.
Since $1-\varepsilon \le (1+\varepsilon )^{-1}\le (1-\varepsilon
)^{-1}\le 1+2\varepsilon $ for $\varepsilon \in [0,\frac12]$, this is
equivalent with the existence of a constant $c_4$ such that
$$
|Cj^{-\alpha }\lambda _{j,+}(P)-1|\le c_4j^{\alpha -\beta },
$$
which is rewritten, with $c_5=C^{-1}c_4$, as
$$
|\lambda _{j,+}(P)-C^{-1}j^\alpha |\le c_5j^{2\alpha -\beta }.\tag5.21
$$

Next we note that the functions $j\to \lambda _{j,+}(P)$
and $t\to N_+(t;P)$ are essentially inverses of one another (in the sense that
$N_+(t;P)$ is a step-function and $j\mapsto\lambda _{j,+}(P)$ should be
filled out at non-integer arguments to have
the reflected graph; both are monotone nondecreasing). To see how one passes
from inequalities for one of them to the other, consider e.g.\ the inequality
$$
\lambda _{j,+}(P)\le C^{-1}j^\alpha +c_5j^{2\alpha -\beta }.
$$
Define $\varphi (j)=C^{-1}j^\alpha +c_5j^{2\alpha -\beta }$.
Let $t=\varphi (j)$ for some $j\in {\Bbb N}$, then 
$$N_+(t;P)\ge
N_+(\lambda _{j,+}(P);P)\ge j.
$$
 Now $t=
C^{-1}j^\alpha +c_5j^{2\alpha -\beta }$ implies $t\le c_6j^\alpha $
(since $2\alpha - \beta <\alpha $) and
$$
(Ct)^{1/\alpha }=(j^\alpha+Cc_5j^{2\alpha -\beta })^{1/\alpha }=
j(1+Cc_5j^{\alpha -\beta })^{1/\alpha }.
$$
Hence
$$
\aligned
j&=(Ct)^{1/\alpha }(1+Cc_5j^{\alpha -\beta })^{-1/\alpha }\ge
(Ct)^{1/\alpha }(1-c_7j^{\alpha -\beta })\\
&\ge (Ct)^{1/\alpha }(1-c_7(c_6^{-1}t)^{(\alpha -\beta )/\alpha })=C^{1/\alpha }t^{1/\alpha }-c_8t^{(1+\alpha -\beta )/\alpha };
\endaligned
$$
for $j$ so large that $Cc_5j^{\alpha -\beta }\le \frac12$; here we
have used the general inequality, valid for $s\in{\Bbb R}$,
$$
1-c_s|x|\le (1+x)^s\le 1+c_s|x|,\text{ for }|x|\le \tfrac12.\tag5.22
$$

This shows that for $t=\varphi (j)$, $j$ sufficiently large,
$$
N_+(t;P)\ge C^{1/\alpha }t^{1/\alpha }-c_8t^{(1+\alpha -\beta )/\alpha },
$$
giving part of the implication from (5.21) to (5.20). The other
needed implications are shown in a similar way. 
 \qed

\enddemo

We shall mainly use the special case where $\alpha =M/n$, $\beta
=(M+\theta )/n$ for some $\theta >0$ and some positive integer $M$, corresponding to $(1+\alpha -\beta )/\alpha =(n-\theta )/M$:

\proclaim{Corollary 5.5} Let $\theta >0$, $C_P>0$. In the setting of Lemma {\rm 5.4}, there exists $c_1>0$ such that 
$$
|\mu _{j,+}(P^{-1})-C_P^{M/n}j^{-M/n }|\le c_1j^{-(M+\theta )/n }\text{ for all
}j\in{\Bbb N},\tag5.23
$$
if and only if there exists $c_2>0$ such that
$$
|N_+(t;P)-C_Pt^{n/M }|\le c_2t^{(n-\theta )/M
}\text{ for all }t>0.\tag5.24
$$

\endproclaim

For the study of the eigenvalues of $A_a$, we note that using the orthogonal decomposition (2.1)
we can write the identity (5.4) in the form
$$
\aligned
A_a^{-1}&=\pr_R\, A_\gamma ^{-1}\pr_R+\pr_R\, A_\gamma
^{-1}\pr_Z+\pr_Z\, A_\gamma ^{-1}\pr_R +\pr_Z\, A_\gamma ^{-1}\pr_Z
+a^{-1}\pr_Z
\\
&=B_1+B_2+S,\text{ with }\\
B_1&=\pr_R\, A_\gamma ^{-1}\pr_R,\\
B_2&=a^{-1}\pr_Z,\\
 S&=\pr_R\,
A_\gamma ^{-1}\pr_Z+\pr_Z\, A_\gamma ^{-1}\pr_R+\pr_Z\, A_\gamma ^{-1}\pr_Z.
\endaligned\tag5.25
$$

For the part $B_1+B_2$, where the two terms act separately in the two
orthogonal subspaces $R$ and $Z$, we see that $B=B_1+B_2$ has the spectrum
$$
\sigma (B_1+B_2)=\sigma (B_1)\cup \sigma (B_2),\tag5.26
$$
consisting of a sequence of positive eigenvalues 
$
\mu _{j,+}(B_1)
$
(since $B_1$ is compact nonnegative), the point 0 (in the essential spectrum) and an
eigenvalue $a^{-1}$ of infinite multiplicity. The essential spectrum
consists of the two points $0$ and $a^{-1}$.
Since $A_a^{-1}$ is a perturbation of $B_1+B_2$ by a compact operator $S$, its essential spectrum again consists of 0 and
$a^{-1}$. As noted earlier, $A_a$ is unbounded above, so it has a sequence of
eigenvalues going to infinity, corresponding to a positive eigenvalue sequence for
$A_a^{-1}$ going to 0.

In the detailed analysis, we shall again take advantage of the calculus of
pseudodifferential boundary operators, 
using some composition rules and an important result
shown in \cite{G84}. The main point is to identify certain terms as
{\it singular Green operators}, which have a better
spectral behavior than the pseudodifferential terms on $\Omega $. We refer to
\cite{G84} for details (introductions to the $\psi $dbo calculus are
also given in \cite{G96} and \cite{G09}).

The following result was shown in \cite{BGW09} Prop.\ 3.5 in the
second-order case:

\proclaim{Proposition 5.6} The orthogonal projection 
$\pr_R$ in $H=L_2(\Omega )$ acts as 
$$
\pr_R=AR_\varrho A =I-\pr_Z,
$$
where $R_\varrho $ is the solution operator for the Dirichlet problem
for $A^2$. Here
$\pr_Z$ is a singular Green operator  on $\Omega $ of order and class $0$.
\endproclaim

\demo{Proof} The proof, formulated in \cite{BGW09} for the nonselfadjoint
second-order case with a spectral parameter, goes over verbatim to the
$2m$-order case, when $\gamma _0,\gamma _1$ are replaced by $\gamma ,\nu $.\qed
\enddemo

In particular, $\pr_R$ and $\pr_Z$ are continuous in $H^s(\Omega )$
for all $s>-\frac 12$.

It follows that all the ingredients in (5.25) are in the $\psi $dbo
calculus:

\proclaim{Proposition 5.7} 

$1^\circ$ The operators
$$
\pr_R\,
A_\gamma ^{-1}\pr_Z,\; \pr_Z\, A_\gamma ^{-1}\pr_R \text{ and }\pr_Z\, A_\gamma ^{-1}\pr_Z,
$$
hence also their sum $S$, cf.\ {\rm (5.25)}, are singular Green
operators on $\Omega $ of order  $-2m$ and class $0$.

$2^\circ$ For any positive integer $N$,
$$
\aligned
&A_\gamma ^{-N}=\pr_R\, A_\gamma
^{-N}\pr_R+S_{1,N},\\
&A_a^{-N}=B_{1,N}+B_{2,N}+S_N,\text{ with }B_{1,N}=\pr_R\, A_\gamma
^{-N}\pr_R,\; B_{2,N}=a^{-N}\pr_Z,
\endaligned
\tag5.27
$$
where $S_{1,N}$ and $S_{N}$ are singular Green operators on $\Omega $ of order  $-2mN$
and class $0$.
\endproclaim 

\demo{Proof} $1^\circ$. It is well-known from the $\psi $dbo calculus that $A_\gamma ^{-1}=A^{(-1)}_++G_\gamma $,
where $A^{(-1)}_+$ is the truncated operator $r^+A^{(-1)}e^+$
and $G_\gamma $ is a singular
Green operator on $\Omega $ of order $-2m$ and class $0$. 
Here $A^{(-1)}$ is a pseudodifferential parametrix
of $A$ extended to ${\Bbb R}^n$, $r^+$ restricts from ${\Bbb R}^n$ to
$\Omega $ and $e^+$ extends by zero on ${\Bbb R}^n\setminus\Omega $. 
Since $\pr_Z$ is a singular Green operator of order and class 0 by
Proposition 5.6, the compositions
with $\pr_Z$ lead to singular Green operators of order $-2m$ and
class 0. Since $\pr_R=I-\pr_Z$, composition with it preserves the order and the property of being a singular
Green operator of class 0. 

$2^\circ$. The statement for the first line of (5.27) has already been shown for 
$N=1$; for general $N$, it follows by similar arguments applied to 
 $A_\gamma
^{-N}$. For the second line of (5.27), we calculate:
$$
\aligned
A_a^{-N}&=(\pr_R\, A_\gamma ^{-1}\pr_R+a^{-1}\pr_Z+S)^N
=(\pr_R\,
A_\gamma ^{-1}\pr_R)^N+a^{-N}\pr_Z+ \text{ s.g.o.s}\\
&=\pr_R\,
A_\gamma ^{-N}\pr_R+a^{-N}\pr_Z+ \text{ s.g.o.s},
\endaligned
$$
by the $\psi $dbo 
rules of calculus, where the s.g.o.s stand for singular Green
operators of class 0 and order $-2mN$. \qed
\enddemo

A main result of \cite{G84} was the following asymptotic estimate of $s$-numbers of
singular Green operators. When $Q$ is a compact operator, its
$s$-numbers are the positive eigenvalues of $|Q|=(Q^*Q)^{1/2}$,
$s_j(Q)=\mu _j(|Q|)$, arranged nonincreasingly and repeated according
to multiplicity.

\proclaim{Theorem 5.8} When $G$ is a singular Green operator on
$\Omega $ of negative order $-M$ and class $0$, then it is compact in
$L_2(\Omega )$ with
$s$-numbers satisfying 
$$
s_j(G)j^{M/(n-1)}\to c(g^0)\text{ for }j\to\infty ,\tag5.28
$$
where $c(g^0)$ is a nonnegative constant defined from the principal
symbol $g^0$ of $G$.
\endproclaim 

The remarkable feature here is that the spectral asymptotics formula
involves the boundary dimension $n-1$ rather than the interior
dimension $n$.

An application to the operators in Proposition 5.7 gives:

\proclaim{Corollary 5.9} The asymptotic property
$$
s_j(G)j^{2mN/(n-1)}\to c(g^0)\text{ for }j\to\infty ,\tag5.29
$$
holds for the singular Green operators $S_N$ and $S_{1,N}$ considered in Proposition
{\rm 5.7}.
\endproclaim

It is seen that $A_a^{-N}$ has several ingredients with different
spectral asymptotics properties. Therefore we need a theorem on
how eigenvalue asymptotics formulas with remainder asymptotics are
perturbed when operators are added together.

This builds on a variant of a result of Ky Fan \cite{F51}.

\proclaim{Lemma 5.10} If $Q$, $B$, and $S$ are bounded
selfadjoint operators whose spectra on ${\Bbb R}_+$ are discrete, 
and $Q=B+S$, then one has for the positive eigenvalues $\mu _{j,+}$, arranged
nonincreasingly and repeated according to multiplicity:
$$
\mu _{j+k-1,+}(B+S)\le \mu  _{j,+}(B)+\mu  _{k,+}(S), 
\tag5.30
$$
for all $j,k$ such that the eigenvalues exist. 

If $S$ has a finite number $K\ge 0$ of positive eigenvalues, then 
$$
\mu _{j+K,+}(B+S)\le \mu  _{j,+}(B), 
\tag5.31
$$
for all $j$ such that the eigenvalues exist. 
\endproclaim

\demo{Proof} The $l$'th positive eigenvalue of $Q$ is characterized by 
$$
\mu _{l,+}(Q)=\min_{u_1,\dots, u_{l-1}\in H}\max\{(Qu,u)\mid \|u\|=1,
\, u\perp u_1,\dots, u_{l-1}\},\tag5.32
$$
as long as this expression is positive; it is reached when the
$u_1,\dots,u_{l-1}$ are an orthogonal system of eigenvectors for the
first $l-1$ positive eigenvalues. Let $x_1,\dots,x_{j-1}$ be
an orthogonal system of eigenvectors for the first $j-1$ positive
eigenvalues of $B$, and let $y_1,\dots,y_{k-1}$ be
an orthogonal system of eigenvectors for the first $k-1$ positive
eigenvalues of $S$.  Then since $Q=B+S$, we have in view of (5.32):
$$
\aligned
\mu _{j+k-1,+}(Q)&\le\max\{(Qu,u)\mid \|u\|=1,
\, u\perp x_1,\dots, x_{j-1}, y_1,\dots,y_{k-1}\}\\
&\le\max\{(Bu,u)\mid \|u\|=1,
\, u\perp x_1,\dots, x_{j-1}\}\\
&\quad +\max\{(Su,u)\mid \|u\|=1,
\, u\perp  y_1,\dots,y_{k-1}\}\\
&=\mu _{j,+}(B)+\mu _{k,+}(S),
\endaligned\tag5.33
$$
showing (5.30). The last statement in case $K=0$ follows from (5.32), since
$(Su,u)\le 0$ then. For $K>0$ it follows from the calculation in (5.33) with $k-1=K$.
\qed
\enddemo

We use this to show, as a variant of \cite{G78} Prop.\ 6.1:

\proclaim{Proposition 5.11}  Let $Q$, $B$, and $S$ be bounded
selfadjoint operators  such
that $Q=B+S$, where the spectrum of  
$B$ in ${\Bbb R}_+$ is discrete, with eigenvalues $\mu _{j,+}(B)\searrow 0$, and $S$ is
compact. Assume
that, with $\beta
>\alpha >0$, $\gamma >\alpha $, and a positive constant $C$,
$$
\align
\mu_{j,+}(B)-Cj^{-\alpha }&\text{ is }O(j^{-\beta })\text{ for }j\to\infty ,\tag5.34\\
s_{j}(S)&\text{ is }O(j^{-\gamma })\text{ for }j\to\infty .\tag5.35
\endalign
$$
Then
$$
\mu_{j,+}(Q)-Cj^{-\alpha }\text{ is }O(j^{-\beta '})\text{ for
}j\to\infty ,\tag5.36 
$$
with 
$$\beta
'=\min\{\beta , \gamma (1+\alpha )/(1+\gamma )\};\tag5.37$$
 here $\beta '\in \,]\alpha ,\beta ]$. 
\endproclaim

\demo{Proof} By hypothesis, $B$ has infinitely many positive
eigenvalues. If $S$ has so too, we proceed as in \cite{G78}, Prop.\ 6.1:
Let $d\in \,]0,1[\,$, to be chosen later. For each $l\in {\Bbb N}$,
let $k=[l^d]+1$ and let $j=l-[l^d]$ in (5.30). Then (5.34)--(5.35)
imply by use of (5.22):
$$
\aligned
\mu _{l,+}(Q)&\le C(l-[l^d])^{-\alpha }+c_2(l-[l^d])^{-\beta
}+c_3([l^d]+1)^{-\gamma }\\
&\le Cl^{-\alpha }(1-[l^d]/l)^{-\alpha }+c_2l^{-\beta
}(1-[l^d]/l)^{-\beta }+c_3l^{-d\gamma }\\
&\le Cl^{-\alpha }+c_2l^{-\beta }+c_4l^{d-\alpha -1} +c_5l^{d-\beta
-1}+c_3 l^{-d\gamma }\\
&\le Cl^{-\alpha }+c_6l^{-\beta '},
\endaligned
$$
where $\beta '=\min\{\beta ,\alpha -d+1,\beta -d+1,d\gamma \}$. Taking
$d=(1+\alpha )/(1+\gamma )$, we have (5.37).

If $S$ has a finite number $K$ of positive eigenvalues, we have if
$K=0$ that 
$$
\mu _{j,+}(Q)\le \mu_{j,+}(B)\le Cj^{-\alpha }+c_1j^{-\beta },\tag5.38
$$
and if $K>0$, for $j\ge K$, by (5.31),
$$
\aligned
\mu _{j,+}(Q)-Cj^{-\alpha }&\le \mu_{j-K,+}(B)
-Cj^{-\alpha  }
\le C(j-K)^{-\alpha } -Cj^{-\alpha  }+c_1(j-K)^{-\beta }\\
&= Cj^{-\alpha }[(1-K/j)^{-\alpha } -1]+c_1j^{-\beta }(1-K/j)^{-\beta
}\\
&\le c_2j^{-\alpha -1} +c_1j^{-\beta }+c_3j^{-\beta -1}\le c_4j^{-\beta ''},
\endaligned\tag5.39
$$
with $\beta ''=\min \{\alpha +1,\beta \}> \beta '$, since $\alpha
+1> \gamma (\alpha +1)/(\gamma +1)$.

This shows the desired upper estimate.
A similar lower estimate is obtained by
noting that Lemma 5.10 applied to $B=Q+(-S)$ gives
$$
\mu _{j,+}(Q)\ge \mu _{j+k-1,+}(B)-\mu _{k,+}(-S).\qquad\square
$$
\enddemo

If needed, one can of course use the finer estimates (5.38) or (5.39)
in appropriate situations. 

The results will first be used to give an
eigenvalue estimate for $\pr_RA_\gamma ^{-N}\pr_R$:

\proclaim{Proposition 5.12} $B_{1,N}=\pr_RA_\gamma ^{-N}\pr_R$ is a
nonnegative compact selfadjoint operator whose positive eigenvalues
satisfy, with $c'_A=c_A^{2mN/n}$, $c_A$ defined by {\rm (5.6)}:
$$
\mu_{j,+}(B_{1,N})-c'_Aj^{-2mN/n }\text{ is }O(j^{-(2mN+\theta _N)/n })\text{
for }j\to\infty, \tag5.40
$$
where $\theta_N =2mN/(2mN+n-1)$.

\endproclaim 

\demo{Proof} Since $\pr_R$ is bounded and $A_\gamma ^{-N}$ is compact,
$B_{1,N}$ is compact. The nonnegativity follows since $A_\gamma ^{-N}\ge
0$ so that
$$
(B_{1,N}u,u)=(\pr_RA_\gamma ^{-N}\pr_Ru,u)=(A_\gamma
^{-N}\pr_Ru,\pr_Ru)\ge 0,
$$
for all $u\in H$.
For the eigenvalue asymptotics, we use the decomposition in the first
line of (5.27),
where $A_\gamma ^{-N}$ has the spectral behavior inferred from (5.7):
$$
\mu _j(A_\gamma ^{-N})-c_A^{2mN/n}j^{-2mN/n}=O(j^{-(2mN+1)/n})\text{
for }j\to\infty ,
$$
 and $S_{1,N}$
has the spectral behavior (5.29), by Corollary 5.9. 
We can then apply Proposition 5.11 with 
$$
\alpha =2mN/n,\quad \beta =(2mN+1)/n,\quad \gamma =2mN/(n-1).\tag5.41
$$
Since 
$$
\tfrac{\gamma (1+\alpha )}{1+\gamma }=\tfrac
{2mN}{n-1}\tfrac{1+2mN/n}{1+2mN/(n-1)}=
\tfrac{2mN+2mN/(2mN+n-1)}n<\beta =\tfrac{2mN+1}n,
$$
we have that
$$
\beta '=(2mN+\theta _N)/n\text{ with } \theta _N=2mN/(2mN+n-1).\qquad\square
$$

\enddemo

Next, we treat the full operator $A_a^{-N}$. The study is easiest to complete when $a<0$.

\proclaim{Theorem 5.13} Consider $A_a^{-N}$; it equals $B+S_N$ with
$B=B_{1,N}+B_{2,N}$ and $S_N$ as in Proposition {\rm
5.7}. Assume that $a<0$. Then when $N$ is odd,  
$$
\mu_{j,+}(A_a^{-N})-c'_Aj^{-2mN/n }\text{ is }O(j^{-(2mN+\theta _N)/n})\text{ for }j\to\infty ,\tag5.42
$$
with $\theta _N=2mN/(2mN+n-1)$, $c_A'= c_A^{2mN/n}$, $c_A$ defined in {\rm (5.6)}.
\endproclaim

\demo{Proof} 
For $B_{1,N}$ we have the asymptotic eigenvalue estimate in
Proposition 5.12.
We add $B_{2,N}$ to $B_{1,N}$, which just adjoins the negative
eigenvalue $a^{-N}$ with infinite 
multiplicity. With $B=B_{1,N}+B_{2,N}$, we now apply
Proposition 5.11 to the sum $A_a^{-N}=Q=B+S_N$, with $\beta =(2mN+\theta _N)/n$. This gives (5.36), with 
$$
\beta '=\min\{\beta , \tfrac {2mN}{n-1 }\tfrac{1+2mN/n}{1+2mN/(n-1)}\}=\beta .\qquad\square
$$  
\enddemo

The cases where $a>0$, or $N$ is even so that $a^N>0$, are handled by
transforming
the problem into one where the
eigenvalue sequence we want to describe runs outside the interval
containing the essential spectrum. 

\proclaim{Theorem 5.14} The conclusion of Theorem {\rm 5.13} holds also
when $N$ is even and when $a>0$.
\endproclaim

\demo{Proof}   It remains to treat the cases where $a^N>0$. Let $b$ be a point in the interval 
$\,]0,a^{-N}[\,$ which is in the
resolvent set of both $B$ and $Q=B+S_N$. Replace $B$ and $Q$ by 
$$
B'=b^2(b-B)^{-1}-b=bB(b-B)^{-1},\quad Q'=b^2(b-Q)^{-1}-b=bQ(b-Q)^{-1}. 
\tag5.43
$$
Then the point $a^{-N}$ in the essential spectrum is moved to
$ba^{-N}(b-a^{-N})^{-1}<0$, whereas the point 0 is preserved, and the sequence of
positive eigenvalues $\mu _{j,+}(B)$ decreasing to 0 in the interval
$\,]0,b[\,$ is turned into the sequence of positive eigenvalues
$$ \mu _{j,+}(B')=b\mu _{j,+}(B)(b-\mu _{j,+}(B))^{-1}\searrow
0.\tag5.44
$$
The operators $B'$ and $Q'$ are of the type treated in Lemma 5.10,
their difference being the compact operator
$$
S_N'=Q'-B'=b^2(b-B-S_N)^{-1}-b^2(b-B)^{-1}=b^2(b-B-S_N)^{-1}S_N(b-B)^{-1}.\tag5.45
$$
Concerning their asymptotic eigenvalue properties, we have that (5.34) implies
$$
\aligned
\mu _{j,+}(B')-Cj^{-\alpha }&=b\mu _{j,+}(B)(b-\mu
_{j,+}(B))^{-1}-Cj^{-\alpha }\\
&=\mu _{j,+}(B)-Cj^{-\alpha }+\mu _{j,+}(B)[b(b-\mu
_{j,+}(B))^{-1}-1]\\
&=\mu _{j,+}(B)-Cj^{-\alpha }+\mu _{j,+}(B)^2(b-\mu
_{j,+}(B))^{-1}\\
&=O(j^{-\beta })+O(j^{-2\alpha }).
\endaligned
\tag5.46
$$
This will be used with $C=c'_A$ and exponents as in (5.40),
$\alpha =2mN/n$ and $\beta =$ \linebreak$(2mN+\theta _N)/n$. Clearly $2\alpha >\beta $, so then
$$
\mu _{j,+}(B')-c'_Aj^{-2mN/n }=O(j^{-(2mN+\theta _N)/n }).
$$
Since $S'_N$ equals $S_N$ composed with  bounded operators, the estimate
(5.35)
implies a similar estimate for $S'_N$. Now Proposition 5.11 can be
applied, with $\alpha $ and $\beta $ as already indicated, and $\gamma =2mN/(n-1)$, showing that the positive
eigenvalues of $Q'$ have the behavior
$$
\mu_{j,+}(Q')-c'_Aj^{-2mN/n }=O(j^{-(2mN+\theta _N)/n})\text{ for }j\to\infty .\tag5.47
$$
Finally this is carried over to the desired behavior of the eigenvalue
sequence $\mu _{j,+}(Q)$ by a calculation similar to (5.46), using
that
$$
Q=bQ'(Q'+b^{-1})^{-1}.\qquad\square
$$

\enddemo

This has the following implications for the counting functions 
for eigenvalues of
$A_a^N$ going to $\infty $:

\proclaim{Theorem 5.15} Let $N$ be a positive integer, and let $r^N>a^N$. The number $N_{+,r^N}(t; A_a^N)$ of eigenvalues of
$A_a^N$ in $[r^N,t]$ behaves
asymptotically as follows:
$$
N_{+,r^N}(t;A_a^N)-c_At^{n/2mN}=O(t^{(n-\theta _N)/2mN})\text{ for }t\to\infty ,\tag 5.48
$$
with $\theta _N=2mN/(2mN+n-1)$, $c_A$ defined by {\rm (5.6)}.

\endproclaim

\demo{Proof} When $a^N<0$, the spectrum of $A_a^N$ is discrete on ${\Bbb
R}_+$, and we can apply Corollary 5.5 directly to (5.42), concluding (5.48) for
$r=0$. A replacement of 0 by some other $r^N>a^N$ only shifts $N_+$ by a fixed
finite number, and does not change the asymptotic property.

Now let $a^N>0$ and take  an $r>|a| $, such that $r^{-N}$ is not in the spectra
of $A_a^{-N}$ and $\pr_RA_\gamma ^{-N}\pr_R$. For this $r$, the number
$N_{+,r^N}(t;A_a^N)$ is the number of eigenvalues of $A_a^N-r^N$ in
$[0,t-r^N]$.

Observe that when we take $b=r^{-N}$ in the proof of Theorem 5.14,
then
$$
Q=A_a^{-N},\quad Q'=r^{-N}A_a^{-N}(r^{-N}-A_a^{-N})^{-1}=(A_a^N-r^N)^{-1}.
$$
For $Q'$ we have the asymptotic estimate (5.47). Then we can apply
Corollary 5.5 to $A_a^N-r^N$ and its inverse $Q'$, concluding that
$$
N_{+,r^N}(t;A_a^N)-c_A(t-r^N)^{n/2mN}=O((t-r^N)^{(n-\theta _N)/2mN})\text{ for }t\to\infty .
$$
This implies (5.48), since $(t-r^N)^s=t^s(1-r^N/t)^s=t^s+O(t^{s-1})$ by (5.22).\qed

\enddemo

We can finally conclude an improved estimate for $A_a$ itself:

\proclaim{Theorem 5.16} Let $r>a$. The number $N_{+,r}(t; A_a)$ of eigenvalues of
$A_a$ in $[r,t]$ behaves
asymptotically as follows, for any $\varepsilon >0$:
$$
N_{+,r}(t;A_a)-c_At^{n/2m}=O(t^{(n-1+\varepsilon )/2m})\text{ for }t\to\infty ,\tag 5.49
$$
with $c_A$ defined in {\rm (5.6)}.
\endproclaim

\demo{Proof} It suffices to consider $r>|a|$.
Since the number of eigenvalues of $A_a$ in $[r,t]$ is the same as the
number of eigenvalues of $A_a^N$ in $[r^N,t^N]$, we conclude from
Theorem 5.15 that
$$
\aligned
N_{+,r}(A_a;t)-c_At^{n/2m}&=N_{+,r^N}(A_a^N;t^N)-c_A(t^N)^{n/2mN}=O((t^N)^{(n-\theta
_N)/2mN})\\
&=O(t^{(n-\theta _N)/2m}).
\endaligned
$$
  Here $N$ can be taken arbitrarily large. Since $\theta
_N=1-(n-1)/(2mN+n-1)\to 1$ for $N\to \infty $, it can for any
  $\varepsilon >0$ be obtained to be
$>1-\varepsilon $, which shows the statement
in the theorem.\qed 
\enddemo

This ends our study of eigenvalue asymptotics for $A_a$, $a\ne 0$. 

\medskip
Actually, some of the above techniques can also be used to  improve the result
of \cite{G83} for $A_0$, so we include this here.

\proclaim{Theorem 5.17} For the discrete eigenvalue sequence of the Krein-von
Neumann extension $A_0$, the number $N_+(t;A_0)$ of eigenvalues in
$\,]0,t]$ satifies, for any $\varepsilon >0$,
$$
N_+(t;A_0)-c_At^{n/2m}=O(t^{(n-1+\varepsilon )/2m})\text{ for }t\to\infty .\tag 5.50
$$
 
\endproclaim

\demo{Proof} We here use some further rules for eigenvalues and
s-numbers, found e.g.\ in Goh\-berg and Krein \cite{GK69}. Denote the positive eigenvalues $\lambda _j(A_0)$,
$j=1,2,\dots$. It is shown in \cite{G83} that their inverses are the
eigenvalues $\mu _j(S_0)$, where  $S_0=R_\varrho ^{1/2}\,A \,R_\varrho
^{1/2}$ as recalled in (5.14); this was used in \cite{G83} to show the
estimate (5.8). Here $R_\varrho ^{\frac12}$ maps $L_2(\Omega )$
bijectively onto $H^{2m}_0(\Omega )$, and the factor $A$ is really
$\Ami$ mapping $H^{2m}_0(\Omega )$ bijectively onto $R=\operatorname{ran}\Ami$,
where one can apply $R_\varrho ^{\frac12}$. They also define mappings
between the spaces intersected with higher-order Sobolev spaces. 

In addition to $S_0$ we shall study iterates of $S_0$.
For $2N$'th powers we can write
$$
S_0^{2N}=(R_\varrho ^{\frac12}A  R_\varrho ^\frac12)^{2N}=
R_\varrho ^{\frac12}(A  R_\varrho )^{2N-1} A  R_\varrho
^\frac12=B_N A  R_\varrho A  \check B_N,
$$
where 
$$
B_N=R_\varrho ^{\frac12}(A  R_\varrho )^{N-1},\quad
\check B_N=(R_\varrho A  )^{N-1} R_\varrho ^{\frac12}.
$$
Here we recognize $AR_\varrho A$ as the projection $\pr_R=I-\pr_Z$,
cf.\ Proposition 5.6. Then
$$
S_0^{2N}=B_N(I-\pr_Z)\check B_N=B_N\check B_N-B_N\pr_Z\check B_N.\tag5.51
$$

The first term is a compact nonnegative operator whose positive eigenvalues
satisfy:
$$
\mu _j(B_N\check B_N)=\mu _j(\check B_N B_N)=\mu _j((R_\varrho
A  )^{N-1} R_\varrho (A  R_\varrho )^{N-1}).
$$
The operator $(R_\varrho
A  )^{N-1} R_\varrho (A  R_\varrho )^{N-1}$ is of the form
$A^{(-2N)}_++{ G}_{2N}$, where 
 $A^{(-2N)}_+$ is the truncation to $\Omega $ of a parametrix   $A^{(-2N)}$ 
of
$A^{2N}$ (as used earlier in the proof of Proposition 5.7), and  ${G}_{2N}$ is a singular Green operator of order $-4mN$ and class
0. 
Then by Corollary 4.5.6
of \cite{G96} we have the asymptotic eigenvalue estimate (in view of Corollary 5.5):
$$
\mu _j(B_N\check B_N)=\mu _j(A^{(-2N)}_++{ G}_{2N})=c'_Aj^{-4mN/n}+O(j^{-(4mN+1-\varepsilon )/n})\text{
for }j\to\infty ,
$$
for any $\varepsilon >0$, with $c'_A=c_A^{4mN/n}$. (It is used here that $A$ is a scalar differential
operator, see the discussion in \cite{G96} Rem.\ 4.5.5 concerning systems.)

For the second term $B_N\pr_Z\check B_N$ we use that there exists a homeomorphism
$$
\Lambda ^{2m}_{-,+}\colon H^{2m+s}(\Omega )\simto H^{s}(\Omega ),\text{ with
inverse }\Lambda ^{-2m}_{-,+}, \text{ any }s\in{\Bbb R},
$$
belonging to the $\psi $dbo calculus, as introduced in \cite{G90} (also
explained in Section 2.5 of \cite{G96}). Then
$$
\aligned
B_N\pr_Z\check B_N&=R_\varrho ^{\frac12}\Lambda ^{2m}_{-,+}\Lambda
^{-2m}_{-,+}(A  R_\varrho )^{N-1}\pr_Z(R_\varrho A 
)^{N-1}\Lambda ^{-2m}_{-,+}\Lambda ^{2m}_{-,+} R_\varrho ^{\frac12}\\
&=
R_\varrho ^{\frac12}\Lambda ^{2m}_{-,+}\widetilde{ G}_{2N}\Lambda
^{2m}_{-,+} R_\varrho ^{\frac12},
\endaligned
$$
where $\widetilde{ G}_{2N}$ is a singular Green operator of order
$-4mN$ and class 0. The operators $R_\varrho ^{\frac12}\Lambda
^{2m}_{-,+}$ and $\Lambda
^{2m}_{-,+} R_\varrho ^{\frac12}$ are bounded in $L_2(\Omega )$. 
Using Theorem 5.8 for $\widetilde{
G}_{2N}$ together with the general rule $s_j(EGF)\le \|E\|
s_j(G)\|F\|$, we find:
$$
s_j(B_N\pr_Z\check B_N)\le C s_j(\widetilde G_{2N})\le C'j^{-4mN/(n-1)}.
$$
 
Now the perturbation result Proposition 5.11 applied to the decomposition in (5.51)
gives (as in the proof of Theorem 5.13):
$$
\mu _j(S_0^{2N})-c_A^{4mN/n}j^{-4mN/n}=O(j^{-(4mN+\theta _{2N} )/n})\text{
for }j\to\infty ;\text{ with }c'_A=c_A^{4mN/n},
$$
with the usual $\theta _{2N}=4mN/(4mN+n-1)$, and hence (as in the
proof of Theorem 5.16)
$$
N_+(t;S_0^{-1})=N_+(t^{2N};S_0^{-2N})=c_At^{n/2m}+O(t^{(n-\theta
_{2N})/2m}), \text{ for }t\to\infty .
$$
Since $\theta _{2N}\to 1$ for $N\to \infty $, and $N$ can be
taken 
arbitrarily large, the assertion of the
theorem  follows. \qed
\enddemo

The validity of the improved estimate (5.50) has been announced by
Mikhailets in \cite{M94}; we have recently been informed that proof
details are in \cite{M06}.

The spectral results in this section  are formulated for a bounded domain $\Omega $ in
${\Bbb R}^n$, but the methods work for general compact manifolds with
boundary, as in \cite{G96}, so the results are valid for such cases too.

\subhead{ Acknowledgement}\endsubhead

The author is grateful to Jan Philip Solovej for useful discussions of
perturbation theorems.



\Refs
\widestnumber\key{[BMNW08]}

\ref\no[AGW11]
\by H. Abels, G. Grubb and I. Wood \paper  Extension theory and  Kre\u\i{}n-type resolvent
  formulas for nonsmooth boundary value
  problems \finalinfo  arXiv:1008.3281, to appear
\endref



\ref\no[AS80]\by A. Alonso and B.  Simon \paper
The Birman-Krein-Vishik theory of selfadjoint extensions of
semibounded operators \jour 
J. Operator Theory \vol 4 \yr1980 \pages 251--270\endref 

\ref\no[AS81]\by A. Alonso and B.  Simon \paper
Addenda to ``The Birman-Krein-Vishik theory of selfadjoint extensions of
semibounded operators'' \jour 
J. Operator Theory \vol 6 \yr1981 \pages 461\endref 


\ref\no[AP04]\by  W. O. Amrein, D. B. Pearson  \paper $M$ operators: a
generalisation of Weyl-Titchmarsh theory\jour  {
J.~Comp.~Appl. Math.}\vol   171\pages 1--26  \yr 2004\endref

\ref\no[A99]\by Yu. M. Arlinskii\paper On functions connected with
sectorial operators and their extensions \jour Integral Equations
Operator Theory \vol 33 \yr 1999\pages 125--152 
\endref

\ref\no[AGMST10] \by M.\ S.\ Ashbaugh, F.\ Gesztesy, M.\ Mitrea, 
R.\ Shterenberg, 
and G.\ Teschl\paper The Krein--von Neumann extension and its connection 
to an abstract buckling problem\jour Math. Nachr.\vol 283 \pages 165--179 \yr2010 
\endref

\ref\no[AGMT10] \by M.\ S.\ Ashbaugh, F.\ Gesztesy, M.\ Mitrea 
and G.\ Teschl\paper Spectral theory for perturbed Krein Laplacians in
nonsmooth domains\jour Adv. Math.\vol 223 \pages 1372--1467 \yr2010 
\endref

 \ref\no[BL07]\by  J.~Behrndt, M.~Langer \paper Boundary value
problems for elliptic partial  differential operators on bounded
domains\jour  { J.~Funct.~Anal.}  \vol  243\pages 536--565 \yr2007 \endref

\ref\no[B56] \by M. S. Birman \paper On the theory of self-adjoint
extensions of positive definite operators\jour 
Mat. Sb. N.S.\vol 38(80) \yr1956 \pages 431--450\finalinfo in Russian \endref

\ref\no[B62]\by M. S. Birman\paper Perturbations of the continuous
spectrum of a singular elliptic operator by varying the boundary and
the boundary conditions
\jour Vestnik Leningrad. Univ. \vol 17 \yr 1962 \pages 22--55
 \transl \nofrills English translation in\book
Spectral theory of differential operators, 
 Amer. Math. Soc. Transl. Ser. 2, 225\publ Amer. Math. Soc.\publaddr
Providence, RI \yr 2008 \pages 19--53  
\endref 
 




\ref\no[BS77] \by M. S. Birman and M. Z.  Solomyak\paper Asymptotic
properties of the spectrum of differential equations
\vol 14 \pages 5--58\jour  Akad. Nauk SSSR Vsesojuz. Inst. Naucn. i
Tehn. Informacii, Moscow\yr 1977 \transl\nofrills English translation
in \jour J. Soviet Math. \vol 12 \yr 1979\pages 247--283 
\endref

\ref\no[BS80] \by   M. S. Birman and M. Z. Solomyak
\paper Asymptotics of the spectrum of variational problems on
solutions of elliptic equations in unbounded domains\jour
  Funkts. Analiz Prilozhen.
  \vol14  \yr1980\pages  27--35\transl\nofrills English translation in
\jour Funct. Anal. Appl. \vol14 \yr1981  \pages267--274
\endref
 

\ref\no[B71]\by 
  L.~Boutet de Monvel  \paper Boundary problems for pseudodifferential
operators\jour  
 {Acta Math.} \vol126\pages  11--51 \yr 1971\endref

\ref\no[BGW09]\by B. M. Brown, G. Grubb, and I. G. Wood \paper $M$-functions for closed
extensions of adjoint pairs of operators with applications to elliptic boundary
problems \jour Math. Nachr. \vol 282\pages 314--347 \yr2009
\endref  

\ref\no[BMNW08]\by 
B.~M.~Brown, M.~Marletta,  S.~Naboko and
 I.~G.~Wood \paper  
Boundary triplets and M-functions for non-selfadjoint operators, with
applications to elliptic PDEs and block operator matrices\jour  
{J. Lond. Math. Soc.}\pages 700--718\vol 77\yr 2008\endref

\comment
\ref\no[B74]\by J. Br\"u{}ning \paper Zur Absch\"a{}tzung der
Spektralfunktion elliptischer Operatoren \jour Math. Z.\vol137
\yr1974\pages 75--85 \endref
\endcomment

\ref\no[BGP06]\by  J.~Br\"{u}ning, V.~Geyler and K.~Pankrashkin \paper
Spectra of self-adjoint extensions and applications to solvable
Schr\"{o}dinger operators\jour Rev. Math. Phys. \vol  20\pages
1-70 \yr 2008\endref 

\ref\no[DM91]\by  V. A.~Derkach and M. M.~Malamud \paper Generalized resolvents and the boundary value problems for Hermitian operators with gaps\jour  { J.~Funct.~Anal.}  \vol  95\pages  1--95 \yr 1991\endref

\ref\no[F51]
\by Ky Fan
\paper Maximum properties and inequalities for the eigenvalues of
completely continuous operators
\jour Proc. Nat. Acad. Sci. USA
\vol 37
\yr 1951
\pages 760-766
\endref


\ref
\no[GM08] \by F. Gesztesy and M. Mitrea \paper Generalized Robin boundary
conditions, Robin-to-Dirichlet maps, and Krein-type resolvent formulas
for Schr\"odinger operators on bounded Lipschitz domains \inbook  Perspectives in Partial Differential Equations, Harmonic Analysis
and Applications: A Volume in Honor of Vladimir G. Maz'ya's 70th Birthday,
Proceedings of
Symposia in Pure Mathematics \eds D. Mitrea and M. Mitrea  
 \vol 79 \publ Amer. Math. Soc.\publaddr Providence,
RI \yr 2008\pages 105--173
\endref

\ref\no[GM11]   
\by F.~Gesztesy and M.~Mitrea
\paper A description of all selfadjoint extensions of the Laplacian
 and Kre\u\i{}n-type
resolvent formulas in nonsmooth domains \jour  J. Analyse Math.
\vol 113 \yr 2011 \pages 53--172
\endref

\ref\no[GK69] \by I. C.  Gohberg and M. G. Krein\book Introduction to the
theory of linear nonselfadjoint operators. Translated from the Russian
by A. Feinstein. Translations of Mathematical Monographs, Vol. 18
\publ American Mathematical Society \publaddr Providence, R.I. \yr
1969 \pages 378  \endref

\ref\no[GM76]\by M. L. Gorbachuk and V. A. Mikhailets\paper Semibounded
selfadjoint extensions of symmetric operators \jour Dokl. Akad. Nauk
SSSR \vol 226\yr 1976\transl\nofrills English translation in \jour
Soviet Math. Doklady\vol 17 \yr1976\pages 185--186
\endref

 \ref\no[GG91]\by  V. I.~Gorbachuk and M. L.~Gorbachuk\book
Boundary value problems for operator differential equations\publ
Kluwer \publaddr Dordrecht \yr 1991
\endref 


\ref\no[G68]\by G. Grubb
\paper A characterization of the non-local boundary value problems
associated with an elliptic operator
\jour Ann\. Scuola Norm\. Sup\. Pisa
\vol22
\yr1968\pages425--513
\endref

\ref\no[G70]\by 
{G.~Grubb} \paper Les probl\`emes aux limites g\'en\'eraux d'un
op\'erateur elliptique, provenant de la th\'eorie variationnelle
 \jour{Bull.~ Sc.~Math.} \vol94\pages 113--157 \yr 1970\endref

\ref\no[G71]\by 
{G.~Grubb} \paper On coerciveness and semiboundedness of general boundary
value problems\jour 
{Israel J. Math.} \vol10\pages 32--95 \yr 1971\endref

\ref\no[G74]\by G. Grubb\paper Properties of normal boundary problems for elliptic
even-order systems\jour Ann\. Scuola Norm\. Sup\. Pisa\vol1{\rm
(ser.IV)}\yr1974\pages1--61
\endref

\ref\no[G78]
\by G. Grubb
\paper Remainder estimates for eigenvalues and kernels of
pseudo-differential elliptic systems
\jour Math. Scand.
\vol43
\yr1978
\pages275--307
\endref

\ref\no[G83]
\by G. Grubb
\paper Spectral asymptotics for the \lq\lq soft" selfadjoint
extension of a symmetric elliptic differential operator
\jour J. Operator
Theory
\vol10
\yr1983
\pages9--20
\endref

\ref 
\key[G84]
\by G. Grubb
\paper Singular Green operators and their spectral asymptotics
\jour Duke Math. J.
\vol 51
\yr 1984
\pages 477--528
\endref

\ref\no[G90] \by G. Grubb\paper Pseudo-differential boundary problems
in $L_p$ spaces\jour Comm. Partial Differential Equations \vol 15 \yr
1990 \pages 289--340\endref

 \ref\no[G96]\by 
{G.~Grubb}\book Functional Calculus of Pseudodifferential
     Boundary Problems
 Pro\-gress in Math.\ vol.\ 65, Second Edition \publ  Birkh\"auser
\publaddr  Boston \yr 1996\endref

\ref\no[G08]
\by G. Grubb
\paper Krein resolvent formulas for elliptic boundary
problems in nonsmooth domains \jour Rend. Sem. Mat. Univ. Pol. Torino
\vol 66 \yr 2008 \pages 13--39 \endref 

\ref\no[G09]\by G. Grubb\book Distributions and operators. Graduate
Texts in Mathematics, 252 \publ Springer \publaddr New York\yr 2009
 \endref
 
\ref\key[G11] \by G. Grubb
\paper
Perturbation of essential spectra of exterior elliptic problems
\jour J. Applicable Analysis \vol 90 \yr 2011 \pages103--123
\endref 

\ref\no[GS95] \by G. Grubb and R. T.  Seeley \paper Weakly parametric
pseudodifferential operators and Atiyah-Patodi-Singer boundary
problems\jour Invent. Math. \vol121 \yr1995 \pages 481--529 \endref

\ref\no[H63]\by 
L.\ H\"ormander \book Linear Partial Differential Operators,
Grundlehren Math. Wiss. vol. 116 \publ Springer Verlag \publaddr
Berlin \yr 1963
\endref


\comment
\ref\no[H68]\by L. H\"o{}rmander\paper The spectral function of an elliptic
operator\jour Acta Math.\vol121\yr1968\pages193--218
\endref
\endcomment

\ref\no[H85]\by L. H\"o{}rmander\book 
 The Analysis of Linear Partial Differential Operators IV,
Fourier Integral Operators, Grundlehren
 Math. Wiss. vol. 275
\publ 
Sprin\-ger Ver\-lag \publaddr Berlin
\yr  1985
\endref

\ref\no[K75]\by  A. N.~Ko\v{c}ube\u{\i} \paper{Extensions of symmetric operators and symmetric binary relations}\jour { Math.~Notes} (1) \vol  17\pages  25--28 \yr 1975\endref


\ref\no[K79]\by V. A. Kozlov \paper Estimation of the remainder in a
formula for the asymptotic behavior of the spectrum of nonsemibounded 
elliptic systems\jour Vestnik
Leningrad. Univ. Mat. Mekh. Astronom. \yr 1979 \pages 112--113 \endref

\ref\no[K83]\by   V. A. Kozlov \paper Remainder estimates in formulas
for the asymptotic behavior of the spectrum for linear operator
pencils\jour  Funktsional. Anal. i Prilozhen. \vol 17 \yr1983 \pages
80--81 
\transl\nofrills English translation in \jour Funct. Analysis
Appl. \vol 17\yr1983 \pages 147--149 \endref

\ref\no[K47] \by M. G. Krein \paper The theory of self-adjoint extensions of
semi-bounded Hermitian transformations and its applications. I
\jour Mat.\ Sbornik \vol 20\pages 431--495 \yr1947\finalinfo  in Russian 
\endref

\ref\no[LM68]\by  J.-L. Lions and E. Magenes \book  Probl\`emes aux
limites non homog\`enes et applications \vol  1 \publ
 \'Editions Dunod \publaddr Paris \yr 1968
\endref

 \ref\no[LS83]\by 
V. E.~Lyantze and  O. G.~Storozh \book Methods of the Theory
of Unbounded
Operators
\publ Naukova Dumka\publaddr  Kiev \yr 1983 \finalinfo  in   Russian
\endref
 

\ref\no[M10] \by M. M. Malamud \paper Spectral theory of elliptic
operators in exterior domains \jour Russian J. Math. Phys.\vol 17 \yr
2010\pages 96--125 \endref

\ref\no[MM02]\by 
{M. M. Malamud} and  {V. I. Mogilevskii} \paper Kre\u\i n
  type formula for canonical resolvents of dual pairs of linear
  relations\jour 
{Methods Funct. Anal. Topology} (4) \vol8\pages 72--100 \yr 2002\endref

\ref\no[M94] \by V. A. Mikhailets \paper Distribution of the
eigenvalues of finite multiplicity of Neumann extensions of an
elliptic operator \jour Differentsial'nye Uravneniya \vol 30 \yr1994
\pages 178--179\transl\nofrills English translation in \jour
Differential Equations \vol 30 \yr1994 \pages 167--168 \endref

\ref\no[M06]
\by V. A. Mikhailets \paper  Discrete spectrum of the extreme nonnegative
extension of the positive elliptic differential operator \inbook
Proceedings of the Ukrainian Mathematical Congress-2001, Section 7,
Nonlinear analysis \publaddr Kyiv \yr 2006 \pages 80--94 
\finalinfo in Russian \endref



\ref
\no[N29] \by J. von Neumann \paper Allgemeine Eigenwerttheorie Hermitescher 
Funktionaloperatoren \jour Math. Ann. \vol 102 \pages 49--131 \yr1929
\endref

\ref\no[P61]\by \AA{}. Pleijel\paper Certain indefinite differential
eigenvalue problems -- the asymptotic distribution of their
eigenfunctions\inbook Part. Diff. Equ. and Continuum Mech. \publ
Wisconsin Press \publaddr Madison\yr1961\pages 19--37
\endref

\ref\no[PR09]   
\by A. Posilicano and L. Raimondi \paper  Krein's resolvent
formula for self-adjoint extensions of symmetric second-order elliptic
differential operators \jour J. Phys. A \vol 42 
  {\rm 015204} \yr 2009 \pages
11 pp \endref


\ref\no[R07] \by V. Ryzhov\paper A general boundary value problem and
its Weyl function \jour Opuscula Math.\vol 27 \yr2007\pages  305--331 
\endref


 \ref\no[V80]\by  
{L. I.~Vainerman}  \paper On extensions of closed operators in Hilbert space\jour 
{Math.~Notes} \vol 28\pages  871--875 \yr 1980\endref

\ref\no[V52]\by
M. I.~Vishik \paper {On general boundary value problems for elliptic
differential operators}\jour 
{Trudy Mosc. Mat. Obsv}  \vol1\pages 187--246 \yr 1952
\transl \nofrills English translation in \jour {Amer. Math. Soc. Transl. (2) } \vol24\pages 107--172 \yr 1963\endref



\endRefs
\enddocument